\documentclass[12pt]{article}
\usepackage[utf8]{inputenc}

\usepackage{pgfplots}
\pgfplotsset{width=7cm,compat=1.16} 
\usepackage{tkz-graph}
\usepackage{tikz}
\usetikzlibrary{decorations.pathreplacing,shapes.misc,calc}
\usetikzlibrary{arrows.meta,decorations.markings,bending}
\usepackage{tkz-euclide}
\usetikzlibrary{arrows.meta}
\tikzset{>=latex} 
\usepackage{xcolor}
\tikzset{
  every point/.style = {circle, inner sep={.75\pgflinewidth}, opacity=1, draw, solid, fill=white},
  point/.style={insert path={node[every point, #1]{}}}, point/.default={},
  colored point/.style = {point={fill=#1}},
  point name/.style = {insert path={coordinate (#1)}},
  inherit/.style = {point/.style={insert path={node[circle, inner sep={.75\pgflinewidth}, draw, fill, #1]{}}}},
  ar3/.style={decoration={markings,mark=at position 0.46 with {\arrow{stealth}},mark=at position 0.97 with {\arrow{stealth}}}, postaction={decorate}},
  ar4/.style={decoration={markings,mark=at position 0.31 with {\arrow{stealth}},mark=at position 0.82 with {\arrow{stealth}}}, postaction={decorate}},
  ar4r/.style={decoration={markings,mark=at position 0.31 with {\arrow[>=stealth]{<}},mark=at position 0.82 with {\arrow[>=stealth]{<}}}, postaction={decorate}},
  ar5/.style={decoration={markings,mark=at position 0.33 with {\arrow{stealth}},mark=at position 0.86 with {\arrow{stealth}}}, postaction={decorate}},
  ->-/.style={decoration={markings,mark=at position 0.15 with {\arrow{stealth}},mark=at position 0.49 with {\arrow{stealth}},mark=at position 0.95 with {\arrow{stealth}}}, postaction={decorate}},
  -<-/.style={decoration={markings,mark=at position 0.18 with {\arrow[>=stealth]{<}},mark=at position 0.5 with {\arrow[>=stealth]{<}},mark=at position 0.88 with {\arrow[>=stealth]{<}}}, postaction={decorate}},
  ar1/.style={decoration={markings,mark=at position 0.13 with {\arrow{stealth}},mark=at position 0.5 with {\arrow{stealth}},mark=at position 0.97 with {\arrow{stealth}}},postaction={decorate}},
  ar2/.style={decoration={markings,mark=at position 0.11 with {\arrow{stealth}},mark=at position 0.51 with {\arrow{stealth}},mark=at position 0.99 with {\arrow{stealth}}},postaction={decorate}}
}

\colorlet{xcol}{blue!60!black}
\colorlet{myred}{red!80!black}
\colorlet{myblue}{blue!80!black}
\colorlet{mygreen}{green!40!black}
\colorlet{mypurple}{red!50!blue!90!black!80}
\colorlet{mydarkred}{myred!80!black}
\colorlet{mydarkblue}{myblue!80!black}
\tikzstyle{xline}=[xcol,thick,smooth]
\tikzstyle{width}=[{Latex[length=5,width=3]}-{Latex[length=5,width=3]},thick]
\tikzstyle{mydashed}=[dash pattern=on 1.7pt off 1.7pt]
\tikzset{
  traj/.style 2 args={xline,postaction={decorate},decoration={markings,
    mark=at position #1 with {\arrow{stealth}},
    mark=at position #2 with {\arrow{stealth}}}
  }
}

\usepackage{graphicx}
\usepackage{hyperref} 
\usepackage{amsmath}
\usepackage{amssymb}
\usepackage{amsthm}
\usepackage{amsfonts} 
\usepackage{bm}
\usepackage{cases}
\usepackage{xcolor}
\usepackage{color,soul}
\usepackage{caption} 
\usepackage{subcaption} 
\usepackage{graphicx}
\usepackage{epstopdf}
\usepackage{enumerate}
\usepackage{ifthen}
\usepackage{mathtools}  
\usepackage{notes2bib}
\usepackage{placeins} 
\usepackage{setspace}
\usepackage{sidecap} 
\usepackage[toc,page]{appendix}
\usepackage{url}
\usepackage{ucs} 

\usepackage{float}

\usepackage{lineno}

\def\epsilon{\varepsilon}
\newcommand{\BF}{\bf\boldmath }

\newcommand{\ie}{{\it{i.e.}}}

\newcommand{\eps}{\varepsilon}

\newcommand{\R}{{\mathbb{R}}}

\newcommand{\bZ}{{\mathbb{Z}}}
\newcommand{\Z}{{\mathbb{Z}}}

\def\cA{{\cal A}}
\def\cC{{\cal C}}
\def\cD{{\cal D}}

\def\cM{{\cal M}}
\def\cN{{\cal N}}
\def\cO{{\cal O}}
\def\cR{{\cal R}}
\def\cP{{\cal P}}
\def\O{{\cal O}}
\def\cQ{{\cal Q}}
\def\cK{{\cal K}}
\def\cS{{\cal S}}

\def\bR{{\mathbb R}}
\def\NW{{\cal P}}

\def\Dt{\cD_{tran}}

\newcommand{\beq}{\begin{linenomath}\begin{equation*}} 
\newcommand{\eeq}{\end{equation*}\end{linenomath}} 
\newcommand{\beqn}{\begin{linenomath}\begin{equation}} 
\newcommand{\eeqn}{\end{equation}\end{linenomath}} 
\usepackage[normalem]{ulem}

\newtheorem*{thmE}{Theorem E}

\newtheorem{theorem}{Theorem}
\newtheorem{proposition}{Proposition}[section]
\newtheorem{definition}[proposition]{Definition}
\newtheorem{lemma}[proposition]{Lemma}
\newtheorem{conjecture}{Conjecture}
\newtheorem{corollary}[proposition]{Corollary}
\newtheorem{example}[proposition]{Example}

\def\Fto{\overset{F}{\geq}}
\def\Lto{\overset{Lyap}{\geq}}
\def\CRto{\overset{CR}{\geq}}
\def\SUMto{\overset{sum}{\geq}}

\def\Pto{\overset{Prol}{\geq}}
\def\leadsto{\Sto}

\def\Fto{\succcurlyeq_F}
\def\NWto{\succcurlyeq_{\cP_F}}

\def\NWtoFG{\succcurlyeq_{\cP_{F_G}}}

\def\NWtoG{\succcurlyeq_{\cP_{\overline F}}}

\def\CRtoI{\succcurlyeq_{\cC_{F,d_1}}}
\def\CRtoII{\succcurlyeq_{\cC_{F,d_2}}}
\def\Sto{\succcurlyeq_S}
\def\CRto{\succcurlyeq_{\cC_F}}
\def\CRKto{\succcurlyeq_{\cC_{F_K}}}
\def\CRGto{\succcurlyeq_{\cC_{F_G}}}
\def\lra{\overset{S}{=}}

\def\Feq{\overset{F}{=}}
\def\Seq{\overset{S}{=}}

\definecolor{dgn}{rgb}{0.2, .6, 0.0}

\newcommand{\black}{\color{black}{}}
\newcommand{\blue}{\color{blue}{}}

\definecolor{orange}{rgb}{0.8, .2, 0.0}

\newcommand{\red}{\color{red}{}}

\definecolor{jim}{rgb}{0,.5,.5}
\definecolor{rob}{rgb}{0,0.4, 0.6}
\def\rob{\color{rob}}

\newcommand{\dist}{{\text{dist}}}
\newcommand{\Dist}{{\text{Dist}}}

\newcommand{\T}{{\mathbb{T}}}

\DeclareMathOperator{\Down}{Down}
\DeclareMathOperator{\Up}{Up}
\DeclareMathOperator{\Node}{Node}
\def\Om{\Omega}

\newcommand{\bydef}{\overset{\mathrm{def}}{=\joinrel=}}

\def\cT{{\cal T}}

\def\sR{\cR}

\usepackage{enumitem}

\newif\ifincludeXX
\includeXXtrue 

\usepackage{lineno}
\newif\ifincludeOLD
\includeOLDtrue 
\includeOLDfalse 
\title{Streams, Graphs\\
and Global Attractors\\ 
of Dynamical Systems\\
on Locally Compact Spaces.}
\author{Roberto De Leo and James A. Yorke}

\begin{document}

\date{\today}
\maketitle

{R. De Leo at
       Department of Mathematics, Howard University, Washington DC 20059,\\
       {roberto.deleo$@$howard.edu}  \\ 
      
      J.A. Yorke at
      Institute for Physical Science and Technology and the Departments of Mathematics and Physics, University of Maryland College Park, MD 20742,\\
      {yorke@umd.edu}
}

\begin{abstract}
    In a recent article, we introduced the concept of streams and graphs of a semiflow.
    An important related concept is the one of semiflow with {\em compact dynamics}, which we defined as a semiflow $F$ with a {\em compact global trapping region}.
    In this follow-up, we restrict to the important case where the phase space $X$ is locally compact and we move the focus on the concept of {\em global attractor}, a maximal compact set that attracts every compact subset of $X$.
    A semiflow $F$ can have many global trapping regions but, if it has a global attractor, this is unique.
    We modify here our original definition and we say that $F$ has compact dynamics if it has a global attractor $G$.
    We show that most of the qualitative properties of $F$ are inherited by the restriction $F_G$ of $F$ to $G$ and that, in case of Conley's chains stream of $F$, the qualitative behavior of $F$ and $F_G$ coincide.
    Moreover, if $F$ is a continuous-time semiflow, then its graph is identical to the graph of its time-1 map.
    Our main result is that, for each semiflow $F$ with compact dynamics over a locally compact space, the graphs of the prolongational relation of $F$ and of every stream of $F$ are connected if the global attractor is connected.
\end{abstract}

\bigskip

\section{Introduction}
We recently introduced~\cite{DLY25} the concept of {\em stream} of a semiflow on a compact metrizable space and showed how to associate to the stream a graph that encodes its main qualitative features.

A stream is a closed and transitive binary relation that establishes ``which point is downstream from which''. 
Given each point $x$, every point on the orbit of $x$ is downstream from $x$ but in general, due to the closure and transitivity, the set of points downstream from a given point is larger than its sole orbit.
Given a stream, it is natural to define ``ponds'' as sets of points that are both upstream and downstream from each other.
Ponds generalize the idea of periodic orbit and we call them {\em nodes} throughout the article because they are the nodes of the graph we associate to the stream.
A pond $M$ is downstream from another pond $N$ if the points of $N$ are downstream from those of $M$.
In this case, we say that there is an edge from $M$ to $N$ in the stream's graph.

The smallest stream was defined in 1964 by Joe Auslander~\cite{Aus63} and the points of its ponds are Auslander's {\em generalized recurrent points}.
The ponds of the smallest stream of an Axiom-A diffeomorphism $f$ are the closed, disjoint transitive invariant subsets of the non-wandering set of $f$ defined by Steven Smale's {\em spectral decomposition}~\cite{Sma67}.

Of course some degree of compactness is needed to grant the existence of ponds. 
In~\cite{DLY25} we introduced the concept of {\em global trapping region} as a forward-invariant set to which all orbits asymptote and such that this convergence is uniform close enough to the region.
Then, we say that a semiflow has {\em compact dynamics} if it has a compact global trapping region.
In this case, there is at least a pond and, correspondingly, the graph is non-empty.

In this article, inspired by a vast literature by the partial differential equations community (e.g. see~\cite{Lad22,Rob01,Hal10,CLR12,Lap18,Lap23}), we introduce a topological concept of {\em global attractor} as a maximal compact set that attracts every compact subset of the system and reformulate all our main concepts and results in terms of global attractors.
One of the advantages of this approach is that, unlike trapping regions, the global attractor of a semiflow, when it exists, is unique.
Moreover, we show that the most important qualitative features of a semiflow with a global attractor are found also in the restriction of the semiflow to its global attractor.
This is why the dynamics of a semiflow with a non-empty global attractor has many features in common with a semiflow over a compact space.
In particular, every such semiflow has at least one node and its graph is non-empty.

For these reasons, we redefine here the concept of compact dynamics by saying that a semiflow has compact dynamics if it has a global attractor. 
This definition update is supported by the following result (see Theorem~\ref{thm: Q <=> G}) we were able to prove: when the phase space is locally compact, a semiflow has a global attractor if and only if it has a compact global trapping region (in the sense of Definition~\ref{def: trapping region})

\smallskip
Within this setting, under the assumption that the phase space $X$ is locally compact and that the semiflow $F$ has compact dynamics, we were able to prove the following main results:
\begin{enumerate}
    \item The following conditions are sufficient for the global attractor of $F$ to be connected: 
    \begin{enumerate}
        \item $F$ has a connected and compact global trapping region (Theorem~\ref{thm: Q conn => G conn};
        \item $F$ is a continuous-time semiflow and has a  path-connected global trapping region (Theorem~\ref{thm: Q path-connected => G connected});
        \item $F$ has a path-connected and locally path-connected global trapping region (Theorem~\ref{thm: pc + lpc => G conn}).
    \end{enumerate}
    \item If the global attractor of $F$ is connected, then the prolongational graph of $F$ (Theorem~\ref{thm:NW}) and the graph of every $F$-stream are connected (Theorem~\ref{thm: connectedness}).
    \item $F$ and its restriction to its global attractor have the same chain-recurrent points, the same chain-recurrent nodes and the same chains graph (Theorem~\ref{thm: CR G}).
    {\em In other words, from the point of view of chains, the qualitative dynamics of $F$ is equivalent to the qualitative dynamics of a semiflow over an invariant compact set.}
    \item Let $F$ be a continuous-time semiflow and denote by $f$ be the time-1 map of $F$ (Theorem~\ref{thm: cont=disc}). 
    Then $F$ and $f$ have the same set of chain-recurrent points, the same chain-recurrent nodes  and the same chains graph.
    {\em In other words, from the point of view of chains, in order to study the full qualitative behavior of a continuous-time semiflow, it is enough to study the behavior of its time-1 map.}
    \item If the Auslander stream has countably many nodes, then it coincides with Conley's chains stream (Theorem~\ref{thm: CF=AF}).
    {\em In particular, the set of generalized recurrent points of $F$ coincides with the set of its chain-recurrent points.}
\end{enumerate}

The article is structured as follows.
In Section~\ref{sec: setting}, we introduce most of the definitions and tools we will use throughout this article. 
In particular, we define global attractors and trapping regions and study the properties that are more relevant to us.
In Section~\ref{sec: connectedness}, we discuss in length about the connectedness of the global attractor depending in the presence of suitable global trapping region for the semiflow.
In Section~\ref{sec: NW}, we study the main properties of the prolongational relation of the non-wandering set.
%
Finally, in Section~\ref{sec: streams}, we study the main properties of streams and, in particular, of chains streams.
\section{Setting, main definitions and basic results}
\label{sec: setting}
\medskip\noindent
{\bf The phase space.} Throughout the article, {\BF$X$} will denote a {\bf metrizable and locally compact} topological space. 
We will usually denote points in $X$ by $x,y,z$ and $d(x,y)$ will denote the distance between $x$ and $y$ for some metric $d$ compatible with the topology of $X$.

\medskip\noindent
{\bf Semi-flows.} The starting point of this work is a discrete-time or continuous-time semi-flow, as defined below.
%
\begin{definition}\label{semi-flow} 
    A {\bf semi-flow} on a topological space $X$ is a continuous map 
    {\BF$F:\T\times X\to X$}, where either $\T=0,1,2,\dots$ (discrete time) or $\T=[0,\infty)$ (continuous time), satisfying the following properties:
    \begin{enumerate}
        \item 
    $F^0(x)=x\text{ for each }x\in X$;
    \item
    $F^{t_1+t_2}(x)=F^{t_2}(F^{t_1}(x))\text{ for each }x\in X\text{ and }t_1,t_2\in\T$.
    \end{enumerate}
    We say that $F$ is a {\bf flow} if, for every $t\geq0$, $F^t$ is invertible.
    In this case, $\T=\bZ$ if the time is discrete and $\T=\R$ if time is continuous; in both cases, we set $F^{-t}=(F^t)^{-1}$.
\end{definition}
\noindent
Notice that the discrete case consists in the iterations of the time-1 map $F^1$.
\begin{definition}[\bf Orbits and limit sets]
    \label{def:orbit}
    Given a semi-flow $F$, we write {\BF $x\Fto y$} if $y=F^t(x)$ for some $t\ge0$ and we say that $y$ is {\BF $F$-downstream from} $x$. We write {\BF $x\Feq y$} if $x\Fto y$ and $y\Fto x$.
    We call {\bf orbit space} of $F$ the set 
    \beq
    \text{\BF$\cO_F$} = \{(x,y): x\Fto y\},
    \eeq
    so that the orbit of any given point $x$ under $F$ is given by
    \beq
    \text{\BF$\cO_F(x)$} = \{y: (x,y)\in\cO_F\}.
    \eeq

    The {\bf limit set} of $F$ is the set
    \beq
    \text{\BF$\Om_F$} = \{(x,y):\text{ there is }t_n\to\infty\text{ as }n\to\infty\text{ such that }F^{t_n}(x)\to y\}.
    \eeq
    The limit set of a point $x$ under $F$ is the set
    \beq
    \text{\BF$\Om_F(x)$}=\{y: (x,y)\in \Om_F\}.
    \eeq
    Similarly, the limit set of a set $A\subset X$ under is the set of points reachable in arbitrarily long times from within $A$:
    \beq
    \text{\BF$\Om_F(A)$} = \{x:\text{ there are }t_n\in\R, a_n\in A\text{ such that }t_n\to\infty, F^{t_n}(a_n)\to x\}.
    \eeq
    
    We say that $x$ is {\bf fixed} for $F$, or that $x$ is a {\bf fixed-point} of $F$, if $\cO_F(x)=\{x\}$; that $x$ is {\bf periodic} if either $x$ is fixed or there is a $y\neq x$ such that $x\Feq y$; that $x$ is {\bf recurrent} if $x\in\Om_F(x)$.
\end{definition}
As illustrated by the example below, the limit set of a set can be strictly larger than the union of the limit sets of its points.
\begin{example}
    Let $F$ be the flow of the ODE $x'=1-x^2$ on $X=[-1,1]$.
    The reader can verify that $\Omega_F(X)=X$.
    On the other side, the limit set of each point in $X$ consists in either the point $-1$ (for $x=-1$) or the point $1$ (otherwise).
\end{example}
The following lemmas will be used several times in the article.
\begin{lemma}
    \label{lemma: O(Q)=cap}
    Let $Q$ be forward-invariant under $F$.
    Then 
    $$
    \Omega_F(Q)=\bigcap_{t\geq0} F^t(Q).
    $$
\end{lemma}
\begin{lemma}
    \label{lemma: unif cont}
    Let $A\subset X$ and assume that $F$ is uniformly continuous on $\T\times A$.
    Then 
    $$\overline{\bigcup_{t\geq0}F^t(A)}=\left(\bigcup_{t\geq0}\overline{F^t(A)}\right)\cup\Omega_F(A).
    $$
\end{lemma}
\begin{proof}
    Let $x\in\overline{\cup_{t\geq0}F^t(A)}$.
    Then there are $a_n\in A$ and $t_n\geq0$ such that $F^{t_n}(a_n)\to x$ as $n\to\infty$.
    Assume first that $t_n$ is bounded. 
    Then, possibly passing to a subsequence and renumbering, there is a $\tau\geq0$ such that $t_n\to\tau$.
    Consider now the sequence $F^\tau(a_n)\in F^\tau(A)$ and notice that
    $$
    d(F^{\tau}(a_n),x)\leq d(F^{\tau}(a_n),F^{t_n}(a_n)) + d(F^{t_n}(a_n),x).
    $$
    Since $F$ is uniformly continuous in $\T\times A$, for every $\eps>0$ we can find a $\delta>0$ such that 
    $|t-t'|+d(a,a')<\delta$ implies
    $d(F^{t}(a),F^{t'}(a'))<\eps/2$.
    Since $t_n\to\tau$ and $F^{t_n}(a_n)\to x$, for every $\eps>0$ we can find an $N>0$ such that $|t_n-\tau|<\delta$ and $d(F^{t_n}(a_n),x)<\eps/2$.
    Hence, $F^\tau(a_n)\to x$ and so $x\in\overline{F^\tau(A)}$.

    Assume now that $t_n$ diverges.
    Then, by definition of limit set, $x\in\Omega_F(A)$, which completes the proof.
\end{proof}
\begin{corollary}
    \label{cor: K}
    Let $K\subset X$ be compact.
    Then 
    $$\overline{\bigcup_{t\geq0}F^t(K)}=\left(\bigcup_{t\geq0}F^t(K)\right)\bigcup\Omega_F(K).
    $$
\end{corollary}
\begin{proof}
    Since $K$ is compact, $F$ is uniformly continuous on $[0,T]\times K$ for every $T>0$.
    Then the same argument of the lemma above proves the claim.
\end{proof}
\begin{definition}
    We call {\bf bitrajectory} of $F$ through $x$ a sequence of points    
    $$b=\{\dots,b_{-1},b_0,b_1,\dots\}$$ 
    such that $F(b_i)=b_{i+1}$ for every $i\in\bZ$ and $b_0=x$.
    We denote by $\alpha(b)$ (resp. $\omega(b)$) the set of limit points of $b$ for $n\to-\infty$ (resp. $n\to\infty$).
\end{definition}
Notice that, if $F$ is a flow, through every point of $X$ passes a unique bitrajectory.
\begin{example}
    \label{ex: bitraj}
    Let $X=\R$ and let $F$ be the flow of any non-zero constant vector field.
    Given any $x$, let $t_n=n$ and $x_n=F^{-t_n}(x)$.
    Then $t_n\to\infty$ and $F^{t_n}(x_n)=x$, so $x\in\Omega_F(X)$.
    Hence, $\Omega_F(X)=X$.
    More generally, one can show in the same way that, if $b$ is any bitrajectory of $F$, then $b\subset \Omega_F(b)$.
\end{example}
The same argument used in the example above proves the following result.
\begin{proposition}
    Let $F$ be a flow.
    Then $\Omega_F(X)=X$.
\end{proposition}
Next example shows that the same phenomenon can happen even in case of semiflows and compact phase spaces.
\begin{example}
    Consider the semiflow given by the logistic map $\ell(x)=4x(1-x)$ on $X=[0,1]$.
    There is a dense orbit in $X$ and so $\Omega_F(X)=X$.
\end{example}
%
\subsection{
Global Attractor, Trapping regions and Compact Dynamics} 
In this article, we continue the study of semi-flows with {\em compact dynamics} we started in~\cite{DLY25}.
In this section we update the definition of compact dynamics we introduced in~\cite{DLY25}, basing it now on the concept of global attractor rather than global trapping region.
\begin{definition}
    Given a set $A$ and a point $x$, we set 
    $$
    d(x,A)=\inf_{a\in A}d(x,a).
    $$
    Given an $\eps>0$ and a set $G\subset X$, we set
    $$
    N_\eps(G)=\{y : d(y,G)<\eps\}.
    $$
    We say that a set $G$ {\bf attracts} a set $K$ under $F$ if, for every $\eps>0$, there exists $T>0$ such that $F^t(K)\subset N_\eps(G)$ for all $t\geq T$.
\end{definition}
\begin{definition}
    The {\bf global attractor} $G\subset X$ of a semiflow $F$ is, when it exists, a maximal invariant compact set of $X$ that attracts each compact set $K\subset X$.
\end{definition}
\begin{proposition}
    Let $G$ be the global attractor of a semiflow.
    Then, for every small enough $\eps>0$, $G$ attracts $N_\eps(G)$.
\end{proposition}
\begin{proof}
    This is a direct consequence of the fact that each compact subset of a locally compact space has a compact neighborhood. 
\end{proof} 
A fundamental property of a global attractor is that it is unique, as shown below.
%
\begin{proposition}
    Assume that $F$ has global attractors $G$ and $G'$. 
    Then $G=G'$. 
\end{proposition}
\begin{proof}
    Since 
    $G$ attracts all compact sets of $X$, for every $\eps>0$ there is $T$ such that $F^t(G')\subset G$ for all $t\geq T$.
    Since $G'$ is invariant, this means that $G'\subset N_\eps(G)$ for every $\eps>0$.
    Hence, $G'\subset G$.
    By the same argument, $G\subset G'$, so that $G=G'$.
\end{proof}
\begin{definition}
    We say that $F$ has {\bf compact dynamics} if it has a global attractor $G$.
\end{definition}
{\BF From now on, throughout the article (unless otherwise specified) $F$ will denote a semiflow with compact dynamics and by $G_F$ its unique global attractor.}
Moreover, unless specified otherwise, all statements in this article hold for both continuous-time and discrete-time semiflows.

\smallskip
The following two proposition illustrate elementary but fundamental properties of global attractors.
\begin{proposition}
    \label{prop: char of G}
    The global attractor $G_F$ contains every compact $F$-inva\-ri\-ant subset of $X$ and is contained in every other set that attracts all compact subsets of $X$.
\end{proposition}
\begin{proof}
    Let $G'$ be a compact $F$-invariant set.
    Then, for each $\eps>0$, there is $\tau>0$ such that $G'=F^t(G')\subset N_\eps(G_F)$ for every $t\geq\tau$, since $G_F$ attracts every compact set.
    Hence, $G'\subset G_F$.
    Assume now that a set $A\subset X$ attracts each compact set. 
    In particular, it attracts $G_F$ and so, by the same argument above,  $G_F\subset A$.
\end{proof}
\begin{lemma}
    \label{lemma: O(A)=A}
    Let $A\subset X$ be a closed $F$-invariant set.
    Then $\Omega_F(A)=A$.
    Assume now that $A$ is compact and let $U\subset X$ be a set such that $U\supset A$ and $A$ attracts $U$.
    Then $\Omega_F(U)=A$.
\end{lemma}
\begin{proof}
    Since $A$ is $F$-invariant, through each point $x\in A$ passes a bitrajectory $b$ (see Example~\ref{ex: bitraj}).
    Since $\Omega_F(A)$ contains all points of its bitrajectories, then $\Omega_F(A)\supset A$.
    Since $A$ is closed, we have also that $\Omega_F(A)\subset A$.
    Hence, $\Omega_F(A)=A$.
    
    Assume now that $A$ is compact and that $U\supset A$. 
    Since $A$ is closed and invariant, $\Omega_F(U)\supset\Omega_F(A)=A$.
    Let $x\in\Omega_F(U)\setminus A$ and set $\eta=d(x,A)$.
    Since $A$ is compact, $\eta>0$.
    Since $A$ attracts $U$, there is a $T>0$ such that $F^t(U)\subset N_{\eta/2}(A)$ for all $t\geq T$.
    Recall that every point in $\Omega_F(U)$ can be arbitrarily approximated in arbitrarily long times with orbits starting from $U$.
    Since all orbits starting within $U$ enter $N_{\eta/2}(G)$ in finite time, assuming that $\Omega_F(U)\setminus A\neq\emptyset$ leads to a contradiction.
    Hence $\Omega_F(U)=A$.
\end{proof}
\begin{proposition}
    \label{prop: omega(K) subset G}
    $\Omega_F(G_F)=G_F$ and $\Omega_F(K)\subset G_F$ for every compact set $K\subset X$.
\end{proposition}
\begin{proof}
    Since $G_F$ is compact and invariant, we know that $\Omega_F(G_F)=G_F$ from Lemma~\ref{lemma: O(A)=A}.
    Let now $K\subset X$ be a compact set and let $x\in\Omega_F(K)$.
    Then there are sequences $t_n\to\infty$ and $x_n\in K$ such that $F^{t_n}(x_n)\to x$.
    Since $G_F$ attracts $K$, then for every $\eps>0$ there is $T_\eps>0$ such that $F^t(K)\subset N_\eps(G_F)$ for all $t\geq T_\eps$.
    Hence, $x\in N_\eps(G_F)$ for every $\eps>0$, namely $x\in G_F$.
    
    Finally, assume that $X$ is locally compact.
    Then, since $G_F$ is compact, $G_F$ has a precompact neighborhood $U$ and, for $\eps>0$ small enough,
    $\overline{N_\eps(G_F)}\subset U$.
    Hence, $\overline{N_\eps(G_F)}$ is compact and so it is attracted by $G_F$.
    Then also $N_\eps(G_F)$ is attracted by $G_F$ and so, by Lemma~\ref{lemma: O(A)=A}, $\Omega_F(N_\eps(G_F))=G_F$. 
\end{proof}
\begin{corollary}
    \label{cor: strong G}
    Let $\eps>0$ be such that $N_\eps(G_F)$ is attracted by $G_F$.
    Then $\Omega_F(N_\eps(G_F))=G_F$.
\end{corollary}
Global attractors are often sets with a highly complicated structure (for instance, they are often not locally connected) and it is in general a hard problem finding out directly their existence. 
In order to at least ascertain their existence, we introduced in~\cite{DLY25} (Definition~2.1.1) the concept of {\em trapping region}.
Unlike global attractors, trapping regions are not unique and often are sets with an elementary structure, such as closed balls.
The definition below updates our previous one in~\cite{DLY25}.
\begin{definition}
    We say that a set $Q\subset X$ {\bf absorbs} a set $K\subset X$ under $F$ if there is a time $T>0$ such that $F^t(K)\subset Q$ for every $t\geq T$.
\end{definition}
\begin{lemma}
    \label{lemma: abs + attr = attr}
    Assume that $Q$ absorbs $U$ and $G$ attracts $Q$.
    Then $G$ attracts $U$.
\end{lemma}
\begin{definition}
    \label{def: trapping region}
    We say that $Q\subset X$ is a {\bf trapping region} for $F$ if $Q$ is forward invariant under $F$ and topologically closed. 
    We say that a trapping region $Q$ is {\bf global} if $Q$ absorbs, under $F$, every compact set $K\subset X$.
    We say that a global trapping region $Q$ is {\bf fat} if $F$ has compact dynamics and $Q$ is a neighborhood of $G_F$.
    We denote by {\BF $\cQ_F$} the set of all global trapping regions of $F$ and by {\BF $\cK_F$} the subset of $\cQ_F$ of the compact global trapping regions.
    Given a trapping region $Q$, we denote by {\BF $F_Q$} the restriction of $F$ to $Q$.
\end{definition}
Notice that $X$ is, trivially, a fat global trapping region for each of its semiflows, so $\cQ_F$ is never empty.
The reader can verify the following elementary properties of $\cK_F$ and $\cQ_F$.
\begin{proposition}
    The sets $\cQ_F$ and $\cK_F$ are invariant under $F^t$ and finite intersections, namely:
    \begin{enumerate}
        \item if $Q\in\cK_F$ then $F^t(Q_1)\in\cK_F$ for all $t\geq0$;
        \item if $Q_1,Q_2\in\cK_F$, then $Q_1\cap Q_2\in\cK_F$;
    \end{enumerate}
    and similarly for $\cQ_F$.
\end{proposition}
%
\smallskip
A fundamental role of trapping regions, as illustrated by the following result, is that one can replace the whole phase space $X$ by any global trapping region of $F$ when studying the global attractor.
{\BF In particular, all results of this article hold, regardless of whether $X$ is locally compact or not, provided $F$ has a locally compact fat global trapping region.}
\begin{proposition}
    \label{prop: G is natural under restrictions}
    Let $F$ have compact dynamics and let $Q\in\cQ_F$.
    Then $F_Q$ has compact dynamics and $G_{F_Q}=G_F$.
\end{proposition}
\begin{proof}
    The set $G_F$ is compact, is invariant under both $F$ and $F_Q$ and, since it attracts all compact sets of $X$, in particular it attracts all compact sets of $Q$.
    Hence, $F_Q$ has a global attractor $G_{F_Q}$ and $G_{F_Q}\supset G_F$.
    On the other side, since $G_F$ is a global attractor for $F$ and $G_{F_Q}$ is compact and invariant under $F$, then $G_{F_Q}\subset G_F$.
    Hence, $G_{F_Q}=G_F$.
\end{proof}
%
%
\begin{example}
    Consider the set $X$ consisting of the disjoint union of a copy of the real line $R$ with the interval $[0,1]$ and let $F$ be a discrete-time semiflow on $X$ such that $F(x)=1$ for each $x\in R$.
    The phase space $X$ is not connected and nor compact.
    These facts, though, play absolutely no role in the dynamics of $F$ since each point, except at most the first point,  of each orbit of $F$ lies in $[0,1]$.
    In this case, $Q=[0,1]$ is a connected and compact global trapping region.
\end{example}
%


%
The following fundamental result shows that every compact global trapping region contains a global attractor and also shows that the definition of ``compact dynamics'' we introduced in~\cite{DLY25} agrees with the one we give here.
\begin{definition}
    Given an $A\subset X$, we say that the set
    $$
    W(A)=\bigcup_{t\geq0}F^t(A)
    $$ 
    is the forward-invariant envelope of $A$.
\end{definition}
\begin{theorem}
    \label{thm: Q <=> G}
    A semiflow $F$ has compact dynamics if and only if $\cK_F\neq\emptyset$.
\end{theorem}
\begin{proof}
    Assume first that $F$ has compact dynamics and
    set 
    $$
    W_\lambda=W(N_\lambda(G_F)),
    \; \lambda>0.
    $$ 
    For $\eps>0$ small enough, $N_\eps(G_F)$ has compact closure and so $G_F$ attracts $N_\eps(G_F)$ and $\Omega_F(N_\eps(G_F))=G_F$ (Corollary~\ref{cor: strong G},).
    By construction, $W_\eps$ contains $G_F$ and is forward-invariant under $F$.
    Moreover, by Lemma~\ref{lemma: unif cont}, for every $\eps>0$ small enough, there is a $T_\eps>0$ such that 
    $$
    \overline{W_\eps} = 
    F^{[0,T_\eps]}\left(\overline{N_\eps(G_F)}\right)\bigcup G_F,
    $$
    so $\overline{W_\eps}$ is compact.
    Hence, $\overline{W_\eps}\in\cK_F$ (notice that, moreover, $\overline{W_\eps}$ is fat).
    
%


    

    Assume now that $F$ has a compact global trapping region $Q$.
    By Lemma~\ref{lemma: O(Q)=cap}, $\Omega_F(Q)=\cap_{t\geq0}F^t(Q)$.
    Since $Q$ is compact, $\Omega_F(Q)$ is non-empty and compact.
    Since $F^t(Q)\subset F^{t'}(Q)$ for $t\geq t'$, then 
    $$
    F^s(\cap_{t\geq0}F^t(Q))=\cap_{t\geq s}F^t(Q)=\cap_{t\geq 0}F^t(Q),
    $$
    namely $\Omega_F(Q)$ is invariant.
    We claim that $\Omega_F(Q)$  attracts $Q$.
    If not, there would be an $\eps>0$ such that $F^t(Q)\not\subset N_\eps(\Omega_F(Q))$ for every $t\geq0$.
    Therefore, we could build a sequence $x_n\in Q$ so that $d(F^n(x_n),\Omega_F(Q))>\eps$ for all $n=1,2,\dots$.
    Since $Q$ is compact, $x_n$ has a subsequence $x_{n_k}$ such that  $F^{n_k}(x_{n_k})\to y$.
    Hence, by definition, $y\in\Omega_F(Q)$.
    By continuity, though, $d(y,\Omega_F(Q))\geq\eps$, that is a contradiction.
    Then, $\Omega_F(Q)$ attracts $Q$ and so $\Omega_F(Q)$ attracts all compact subsets of $X$.
    Moreover, by construction, there is no invariant set larger than $\Omega_F(Q)$. 
    Hence, $\Omega_F(Q)$ is the global attractor of $F$ and so $F$ has compact dynamics.
\end{proof}
\begin{proposition}
    \label{prop: F has a fat compact global trapping region}
    For each $\eps>0$, there is a fat $Q\in\cK_F$ with $Q\subset N_\eps(G_F)$.
\end{proposition}
\begin{proof}
    The argument used to prove Theorem~\ref{thm: Q <=> G} shows that, for $\delta>0$ small enough, $\overline{W(N_\delta(G_F))}$ is a fat compact global trapping region for $F$. 
    Now let $\eta\in(0,\delta)$ and suppose that, for every $\rho>0$, $W(N_\rho(G_F))\not\subset N_\eta(G_F)$.
    Then there are sequences $x_n\in N_{1/n}(G_F)$ and $t_n\geq0$ such that $d(F^{t_n}(x_n),G_F)>\eta$ for all $n=1,2,\dots$.
    We can assume without loss of generality that $x_n\to\bar x\in G_F$.
    If $t_n$ is bounded, then we can assume without loss of generality that $t_n\to\bar t$ and so $F^{t_n}(x_n)\to F^{\bar t}(\bar x)\in G_F$.
    By continuity, though, $d(F^{\bar t}(\bar x),G_F)\geq\eta$, which contradicts the fact that $F^{\bar t}(\bar x)\in G_F$.
    If $t_n$ is unbounded, then nevertheless $F^{t_n}(x_n)\to y\in\Omega_F(N_\rho(G_F))=G_F$, leading to the same contradiction.
    Hence, for every $\eta>0$ there is a $\rho>0$ such that $W(N_\rho(G_F))\subset N_\eta(G_F)$.
    Then $\overline{W(N_\rho(G_F))}\subset \overline{N_\eta(G_F)}\subset N_\eps(G_F)$ for every $\eps>\eta$, which proves the claim.
\end{proof}
\begin{corollary}
    \label{cor: G attracts Q}
    Let $Q$ be a compact global trapping region of $F$.
    Then
    $$
    G_F=\bigcap_{t\geq0}F^t(Q)
    $$
    and $G_F$ attracts $Q$.
\end{corollary}
\begin{proposition}
    $$
    G_F=\bigcap_{Q\in\cQ_F}Q=\bigcap_{Q\in\cK_F}Q.
    $$
\end{proposition}
\begin{proof}
    Let $Q\in\cQ_F$.
    Recall that $Q$ absorbs every compact subset of $X$.
    Since $G_F$ is compact and invariant under $F$, the only possibility is that $G_F\subset Q$, so that $G_F\subset\cap_{Q\in\cQ_F}Q$.

    By Theorem~\ref{thm: Q <=> G}, $F$ has a compact global trapping region $Q_0$ and $\Omega_F(Q_0)=G_F$.
    Notice that, if $Q_0$ is a compact global trapping region, then also $F^t(Q_0)$ is for every $t\geq0$, so that 
    $$
    G_F=\Omega_F(Q_0)=\bigcap_{t\geq0}F^t(Q_0)\supset \bigcap_{Q\in\cQ_F}Q.
    $$
    Hence, $G_F=\cap_{Q\in\cQ_F}Q$.
\end{proof}
Next result improves, respectively, Proposition~2.1.2 and~2.1.3 in~\cite{DLY25}.
%
\begin{proposition}
    \label{prop:recurrent pt}
    The following hold:
    \begin{enumerate}
        \item for each $x$, $\Om_F(x)$ is a non-empty subset of $G_F$;
        \item for every $x$, there is a $F$-recurrent point in $\Omega_F(x)$.
    \end{enumerate}
\end{proposition}
\section{Connectedness of the global attractor}
\label{sec: connectedness}
The connectedness of the global attractor is very important to us since our main result is that it implies the connectedness of the graph of the prolongational relation (Theorem~\ref{thm:NW}) and of every stream (Theorem~\ref{thm: connectedness}) of a semiflow with compact dynamics on a locally compacted space.

In this section we present several conditions that grant the connectedness of the global attractor.
Our first result depends solely on the existence of a suitable trapping region.
\begin{theorem}
    \label{thm: Q conn => G conn}
    If $F$ has a connected compact global trapping region, then $G_F$ is connected.
\end{theorem}
\begin{proof}
    Suppose that $G_F=G_1\cup G_2$ with $G_1$ and $G_2$ compact and mutually disjoint and let $U_1$ and $U_2$ two disjoint neighborhoods of, respectively, $G_1$ and $G_2$.
    Then $U_1\cup U_2$ is a neighborhood of $G$ and so there is an $\eps>0$ such that $N_{\eps}(G_F)\subset U$.
    
    Since, by Corollary~\ref{cor: G attracts Q}, $G_F$ attracts $Q$, there is a $\tau>0$ such that $F^t(Q)\subset N_\eps(G_F)$ for all $t\geq\tau$.
    Let $\bar t$ such that $F^{\bar t}(Q)\subset U$.
    Since $F^{\bar t}(Q)$ is connected, then either $F^{\bar t}(Q)\subset U_1$ or $F^{\bar t}(Q)\subset U_2$.
    Hence, either $G_F\subset U_1$ or $G_F\subset U_2$, contradicting the assumption that $G_1\subset U_1$ and $G_2\subset U_2$.
\end{proof}
The reminder of the section is inspired by an article by M.~Gobbino and M.~Sardella~\cite{GobS97} on the connectedness of another type of global attractors, introduced and widely used in the theory of PDEs (hence, in a not locally compact setting), defined as maximal compact sets that attract all bounded subsets of $X$ (notice that this definition is not topological).
Clearly, those global attractors are also global attractors with respect to our definition.
%
\begin{lemma} 
    \label{prop: invariant connected components}
    Let $F$ be a continuous-time semiflow with compact dynamics.
    Then each connected component of $G_F$ is invariant under $F$.
\end{lemma}
\begin{proof}
    Let $G_F=G_1\cup G_2$ with $G_1,G_2\subset G_F$ compact and mutually disjoint.
    Let $x\in G_1$.
    Since $F^{[0,\infty)}(x)$ is connected and has a point in $G_1$, it must be a subset of $G_1$.
    Hence $G_1$ is forward invariant.
    Let now $y\in G_F$ be such that $F^\tau(y)=x$.
    Since $F^{[0,\tau]}(y)$ is connected and $F^\tau(y)\in G_1$, then $F^{[0,\tau]}(y)\subset G_1$.
    Hence $G_1$ is also backward invariant and so is invariant.
    The same holds for $G_2$.
\end{proof}
We point out that this property does not hold for discrete-time semiflows.
Consider, for instance, the time-1 map $f=F^1$ of the flow $F$ of the ODE $x'=(1-x^2)$ and set $X$ to be the union of the fixed points $0$ and $1$ with a single two-sided trajectory $(\dots,x_{-1},x_0,x_1,\dots)$ of $f$.
In particular, $x_n<x_{n+1}$ for every $n\in\Z$ and $\lim_{n\to\pm\infty}=\pm1$.
The global attractor $G_F$ is $X$ itself and each point of the two-sided trajectory is a connected component of $G_f$.
Yet, no such connected component is invariant since, as pointed out above, $f(x_n)=x_{n+1}$.
\begin{theorem}
    \label{thm: Q path-connected => G connected}
    Assume that a continuous-time semiflow $F$ with compact dynamics has a path-connected global trapping region $Q$.
    Then $G_F$ is connected.
\end{theorem}
Notice that here we are not assuming $Q$ to be fat or compact.
In particular, $G_F$ might not attract $Q$. 
\begin{proof}
    Assume that $G_F$ is not connected. 
    Then there must be two mutually disjoint compact sets $G_1,G_2\subset G_F\subset Q$ such that $G_F=G_1\cup G_2$.
    As argued in Lemma~\ref{prop: invariant connected components}, 
    $G_1$ and $G_2$ are invariant.
    Moreover, since they are compact, there is a $\eps>0$ small enough that $N_\eps(G_1)$ and $N_\eps(G_2)$ are disjoint.
    Let $x_1\in G_1$ and $x_2\in G_2$.
    Since $Q$ is path-connected, there is a continuous path $\gamma:[0,1]\to Q$   from $x_1$ to $x_2$.
    Since $\gamma([0,1])\subset Q$ is compact and $G_F$ attracts every compact subset of $Q$, there is a $\tau>0$ such that $F^t(\gamma([0,1]))\subset N_\eps(G_F)$ for all $t\geq\tau$.
    Since $\gamma([0,1])$ is connected and $F^t$ is continuous, then  $F^t(\gamma([0,1]))$ is connected as well and so either $F^t(\gamma([0,1]))\subset N_\eps(G_1)$ or $F^t(\gamma([0,1]))\subset N_\eps(G_2)$.
    Neither of those two possibilities can arise, though, since $\gamma(0)=x_1\in G_1$ and $\gamma(1)=x_2\in G_2$.
    Since the assumption of $G_F$ being not connected leads to a contradiction, $G_F$ must be connected.    
\end{proof}
In case of discrete-time semiflows, the following property, adapted from a  result by Gobbino and Sardella in~\cite{GobS97}, plays a relevant role.

\begin{proposition}
    \label{prop: finitely many components}
    Assume that there is $Q\in\cQ_F$ that is path-connected global trapping region.
    Then either $G_F$ is connected or it has infinitely many connected components.   
\end{proposition}
\begin{proof}
    By Proposition~\ref{prop: G is natural under restrictions}, we can restrict without loss of generality $F$ to $Q$.
    We assume, by contradiction, that $G_F$ has a finite number $m\geq2$ of connected components $G_1,\dots,G_m$.
    Since $m$ is finite, there is an $\eps>0$ such that the sets $N_\eps(G_i)=\{y\in Q:d(y,G_F)<\eps\}$, $i=1,\dots,m$, are all disjoint subsets of $Q$.
    Moreover, since $G_F$ is invariant under $F$, there is a permutation $\sigma$ of $\{1,\dots,m\}$ such that $F^1(G_i)=G_{\sigma(i)}$.
    In particular, then, $F^{k\cdot m!}(G_i)=G_i$ for every $i=1,\dots,m$ and every $k\geq0$.

    Let now $x_1\in G_1$ and $x_2\in G_2$ and denote by $\gamma:[0,1]\to Q$ a continuous path from $x_1$ to $x_2$.
    Since $\gamma([0,1])$ is compact and $G$ attracts every compact set of $Q$, for every $\eps>0$ there is a $n_\eps\geq0$ such that $F^{n}(\gamma([0,1]))\subset N_\eps(G_F)$ for $n\geq n_\eps$.
    This leads to the following contradiction: $\gamma([0,1])$ is connected and so, for $k$ so large that $k\cdot m!>n_\eps$, the set $F^{k\cdot m!}(\gamma([0,1]))$ must belong to the $\eps$-neighborhood of a single $G_i$; on the other side, $F^{k\cdot m!}(x_1)\in G_1$ and $F^{k\cdot m!}(x_2)\in G_2$.    
\end{proof}
In order to prove our result on the connectedness of the global attractor of a discrete-time semiflow, we need the following lemma.
\begin{lemma}[Gobbino and Sardella, 1997~\cite{GobS97}]
    \label{lemma: GS}
    Let $U\subset X$ be a set with $m$ connected components such that:
    \begin{enumerate}
        \item $\Omega_F(U)$ is compact and attracts $U$;
        \item $\Omega_F(U)\subset U$.
    \end{enumerate}
    Then $\Omega_F(U)$ has at most $m$ connected components.
\end{lemma}
\begin{theorem}
    \label{thm: pc + lpc => G conn}
    Assume that $F$ has a path-connected and locally path-connected global trapping region.
    Then $G_F$ is connected.
\end{theorem} 
\begin{proof}
    Let $Q\in\cQ_F$ be path-connected and locally path-connected.
    By Proposition~\ref{prop: G is natural under restrictions}, we can restrict $F$ to $Q$.
    Since $G_F$ is compact, for every $\eps>0$ there are fintiely many points $x_1,\dots,x_k\in G_F$ such that $G_F\subset\cup_{i=1}^k N_\eps(x_i)$, where $N_\eps(x_i)=\{y\in Q: d(y,x_i)<\eps\}$.
    Since $Q$ is locally path-connected, for $\eps>0$ small enough each of the $N_\eps(x_i)$ is connected and so the set $U=\cup_{i=1}^k N_\eps(x_i)$ has finitely many connected components.
    By construction, $U$ is a compact neighborhood of $G_F$.
    Hence, by Lemma~\ref{lemma: O(A)=A}, $\Omega_F(U)=G_F\subset U$.
    Then, by Lemma~\ref{lemma: GS}, $G_F$ has finitely many connected components.
    Finally, by Proposition~\ref{prop: finitely many components}, $G_F$ consists in a single connected component. 
\end{proof}
We summarize our results on the connectedness of the global attractor in the following corollary.
\begin{corollary}
    Each of the conditions below is sufficient for the connectedness of $G_F$:
    \begin{enumerate}
        \item $F$ is a continuous-time semiflow and has a path-connected global trapping region;
        \item $F$ has a path-connected and locally path-connected global trapping region;
        \item $F$ has a connected  compact global trapping region.
    \end{enumerate}
\end{corollary}
%

\subsection{A compact dynamics semiflow on a connected space with a not connected global attractor}
The following example by M. Gobbino and M. Sardella~\cite{GobS97} shows the non-triviality of the results above on the connectedness of the global attractor.
The example shows a 
discrete-time semiflow on a connected (but not locally connected) phase space $X$ whose global attractor is not connected.
\begin{figure}
    \centering
    \includegraphics[width=\linewidth]{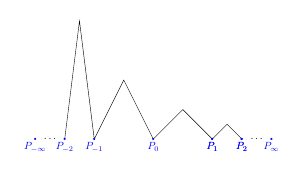}
    \caption{Caption}
    \label{fig:enter-label}
\end{figure}

Let $P_n$, $n\in\Z$, be a sequence of points in the real line for which there exists two points $P_{\pm\infty}$ such that:
$$
P_n>P_{n-1},\;n\in\Z;\;
\lim_{n\to\pm\infty} P_n=P_{\pm\infty}.
$$
Consider first the case $X=P_\infty\cup P_{-\infty}\cup_{n\in\Z}P_n$ with the semiflow $F$ defined by 
$F(P_{\pm\infty})=P_{\pm\infty}$, $F(P_n)=P_{n+1}$.
Then the global attractor of $F$ is the whole phase space $X$, in particular it has infinitely many connected components.
The lack of connectedness of the global attractor is not particularly surprising because $X$ itself has infinitely many connected components.
As Gobbini and Sardella showed, though, an elementary modification of $X$ and $F$ results in a connected phase space with exactly the same global attractor.

Let us embed $X$ in $\R^2$ as a subset of the $x$ axis.
We call $P_n$ the image in the plane of the points $P_n$ since there will be no ambiguity.
Denote by $T_n$ the isosceles triangle of height $2^{-n}$ with basis the segment with endpoints $P_n$ and $P_{n+1}$ and third vertex in the half-plane above the $y$ axis.
Denote by $X_n$ the union of its two sides of $T_n$ of equal length.
Finally, denote by $X'$ the union of $X$ with the sets $X_n$, $n\in\Z$ and  define $F'$ as the map that coincides with $F$ on $X$ and sends piecewise-linearly $X_n$ to $X_{n+1}$.
In particular, $F'$ acts on the second coordinate of every point of $X'$ as the division by 2.
The reader can verify that $F'$ is a discrete-time flow on $X'$.

We claim that the global attractor of $F'$ is $X$.
Indeed, let $Q$ be the intersection of $X'$ with any rectangle of finite height containing $X$.
Then $Q$ is a global trapping region of $F'$, since each $X_n$ is absorbed by $Q$ in finite time.
Moreover, $Q$ is compact and so the global attractor equals $\cap_{n\geq0}F^n(Q)$.
The reader can easily verify that the points in $X$ are indeed the only ones that belong to $F^n(Q)$ for every $n\geq0$.
\section{The prolongational relation $\NW_F$}
\label{sec: NW}
\begin{definition}
    We denote by {\BF$\NW_F$} the relation $\overline{\cO_F}$ and by {\BF$\NWto$} the corresponding symbol. 
    We call this relation the {\bf prolongational relation}.
\end{definition} 
\begin{example}
    Let $F$ be the flow of the ODE $x'=-\sin(\pi x)$ on $X=[0,1]$.
    The orbit space $\cO_F\subset X^2$ is the triangle with vertices $(0,0)$, $(1,0)$, $(1,1)$ minus the boundary points $(1,x)$, $x\in[0,1)$.
    Notice that $1\Fto x$ if and only if $x=1$, since $1$ is fixed.
    The prolongational relation $\overline{\cO_F}$ is the closed triangle with the vertices mentioned above.
    Hence, $1\NWto x$ for every $x\in[0,1]$.
\end{example}
\begin{definition}
    Given points $x,y$ and an $\eps>0$, a {\BF $(F,\eps)$-link} from $x$ to $y$ of length $n+1$ is a finite orbit segment $(\zeta,F(\zeta),\dots,F^n(\zeta))$ such that $d(x,\zeta)<\eps$ and $d(y,F^n(\zeta))<\eps$.
    Given a $(F,\eps)$-link $(\zeta_1,F(\zeta_1),\dots,F^{n_1}(\zeta_1))$ from $x$ to $y$ and a second $(F,\eps)$-link $(\zeta_2,F(\zeta_2),\dots,F^{n_2}(\zeta_2))$ from $y$ to $z$, we say that the two $(F,\eps)$-links are {\em linkable} if $F^{n_1}(\zeta_1)=\zeta_2$.
\end{definition}
Given the two linkable $(F,\eps)$-links above, the sequence 
$$
(\zeta_1,\dots,F^{n_1}(\zeta_1),F(\zeta_2),\dots,F^{n_2}(\zeta_2))
$$ 
is a $(F,\eps)$-link from $x$ to $z$.
\begin{proposition}
    $x\NWto y$ if and only if, for every $\eps>0$, there is an $(F,\eps)$-link from $x$ to $y$.
\end{proposition}
\begin{definition}
    \label{def:NW}
    We say that a point $x$ is {\bf non-wandering} for $F$ if, for every $\eps>0$, there is an $(F,\eps)$-link from $x$ to itself.
    We say that $x,y$ are {\BF $\NW_F$-equivalent} if $x\NWto y$, $y\NWto x$ and, for every $\eps>0$, there is a pair of linkable $(F,\eps)$-links from $x$ to $y$ and from $y$ to $x$.
    We denote by {\BF$NW_F$} the set of all non-wandering points of $F$.
    We say that a set $M\subset NW_F$ is $\NW_F$-equivalent if all points of $M$ are $\NW_F$-equivalent to each other.
    We call {\bf nodes} of $NW_F$ the maximal $\NW_F$-equivalent subsets of $NW_F$.
\end{definition}
The $\NW_F$-equivalence induces a decomposition of $NW_F$ as follows.
\begin{definition}
      We call {\bf node} of $NW_F$ each maximal $\NW_F$-equivalent subset of $NW_F$.  
\end{definition}
Below we recall some fundamental dynamical property of the non-wandering set and its nodes from~\cite{DLY25}.
\begin{proposition}
    \label{prop: Om(x)}
    The following hold:
    \begin{enumerate}
        \item For every bitrajectory $b$, $\alpha(b)$ and $\omega(b)$ are $\NW_F$-equivalent sets (not necessarily $\NW_F$-equivalent to each other).
        \item For every $x$, the set $\Omega_F(x)$ is $\NW_F$-equivalent.
        \item If $N$ is a node of $NW_F$ and $\Omega(x)\cap N\neq\emptyset$, then $\Omega(x)\subset N$.
        \item If $x\in NW_F$, the set $\cO_F(x)\cup\Omega_F(x)$ is $\NW_F$-equivalent.
        \item If $x\in NW_F$ belongs to a node $N$, then $\cO_F(x)\cup\Omega_F(x)\subset N$.
        \item $NW_F$ and each of its nodes are closed and forward-invariant under $F$.
    \end{enumerate}
\end{proposition}
%

%

%
Under certain conditions, the non-wandering set and all of its nodes are invariant. 
Below we present two general conditions under which this holds.
\begin{lemma}
    \label{lemma: NW_F=NW_G}
    Let $F$ be a flow and denote by $\overline F$ the inverse flow, namely the flow $\overline F:\T\times X\to X$ defined by ${\overline F}^t(x)=F^{-t}(x)$.
    Then:
    \begin{enumerate}
        \item $x\NWto y$ if and only if $y\NWtoG x$.
        \item $NW_F=NW_{\overline F}$.
        \item $x,y$ are $\NW_F$-equivalent if and only if they are $\NW_{\overline F}$-equivalent.
        \item $N$ is a node of $NW_F$ if and only if it is a node of $NW_{\overline F}$.
    \end{enumerate}
\end{lemma}
\begin{proof}
    (1) Assume that $x\NWto y$.
    Then, for every $\eps$, there is a $(F,\eps)$-link $$(z,F(z),\dots,F^n(z))$$ from $x$ to $y$.
    Let now $w=F^n(z)$.
    Then ${\overline F}^k(w)=F^{n-k}(z)$ for all $k\in\T$ (in particular, ${\overline F}^n(w)=z$) and so the orbit segment $$(w,{\overline F}(w),\dots,{\overline F}^n(w))$$ is a $({\overline F},\eps)$-link from $y$ to $x$.
    Hence, $y\NWtoG x$.
    The same argument applied to the case  $y\NWtoG x$ proves the claim. 

    (2) This case can be proved using the same argument used in case (1).

    (3) By case (1), if $x\NWto y$ and $y\NWto x$ then also $y\NWtoG x$ and $x\NWtoG y$.
    Moreover, if two $(F,\eps)$-links are linkable, then the corresponding $(G,\eps)$-links defined in case (1) are linkable as well.
    By Definition~\ref{def:NW}, this means that if $x$ is $\NW_F$-equivalent to $y$ then $x$ is $\NW_G$-equivalent to $y$ and viceversa.

    (4)
    This is an immediate consequence of (3).
\end{proof}
\begin{proposition}
    Let $F$ be a flow.
    Then $NW_F$ and all of its nodes are invariant under $F$.    
\end{proposition}
\begin{proof}
    Let $N$ be a node of $NW_F$.
    We know from Proposition~\ref{prop: Om(x)} that each node of $NW_F$ is forward-invariant under $F$ and we know from Lemma~\ref{lemma: NW_F=NW_G} that each node of $NW_F$ is also forward-invariant under the inverse flow $G^t=(F^t)^{-1}$.
    Hence, if $x\in N$, the whole (unique) bitrajectory $\cO_F(x)\cup\cO_G(x)$ passing through $x$ lies in $N$.
    Hence, $F(N)=N$, i.e. $N$ is invariant under $F$.
    Being the union of invariant sets, $NW_F$ is invariant as well.
\end{proof}
\begin{proposition}
    If $F^t$ is an open map for every $t\in\T$, then $NW_F$ and all of its nodes are invariant under $F$.    
\end{proposition}
\begin{proof}
    Let $x\in NW_F$.
    If $x$ is periodic, the claim is trivial, so we can assume that $F^k(x)\neq x$ for every $k\geq0$.
    Let $\eps_n>0$, $n=1,2,\dots$, be a sequence such that $\eps_n\to0$.
    Then, for every $n$, there exist a $(F,\eps_n)$-link $(z_n,F(z_n),\dots,F^{k_n}(z_n))$ from $x$ to itself, i.e. $d(x,z_n)<\eps_n$ and $d(x,F^{n_k}(z_n))<\eps_n$.
    Since $x$ is not periodic, the length of the links diverges as $\eps_n\to0$, namely $k_n\to\infty$.
    Since $F$ has compact dynamics, by Theorem~\ref{thm: Q <=> G}, $F$ has a compact global trapping region $Q$.
    Then, since $X$ is locally compact, for $\eps>0$ small enough the set $N_\eps(Q)$ is compact as well and so there is some $T>0$ such that $F^t(x)\in Q$ for each $x\in N_\eps(Q)$ and $t\geq T$.
    Notice that, for almost all $n$, $z_n\in N_\eps(Q)$ and $k_n-1\geq T$.
    Hence, for almost all $n$, $F^{k_n-1}(z_n)\in Q$.
    Let $y_n=F^{k_n-1}(z_n)$.
    Since $Q$ is compact, we can assume (possibly passing to a subsequence) that $y_n\to\bar y\in Q$.
    Since $F(y_n)=F^{k_n}(z_n)\to x$, by continuity we get that $F(\bar y)=x$.
    In particular, $\bar y\NWto x$.
    Notice also that the orbit segment $(z_n,F(z_n),\dots,F^{k_n-1}(z_n))$ is a $(F,\eps_n)$-link from $x$ to $y_n$, so that $x\NWto y_n$.
    Hence, since $NW_F$ is closed, $x\NWto\bar y$.

    Finally, since $F^1$ is an open map, for each $n=1,2,\dots$ we can find a $w_n$ so that $F^1(w_n)=z_n$ and $w_n\to\bar y$.
    Hence, the orbit segment
    $$
    (w_n,z_n,F(z_n),\dots,y_n)
    $$
    is a $(F,\eps_n)$-link from $\bar y$ to itself, namely $\bar y\in NW_F$.
    Moreover, the $(F,\eps_n)$-link $(w_n,z_n)$ from $\bar y$ to $x$
    is linkable to the $(F,\eps_n)$-link $(z_n,\dots,y_n)$ from $x$ to $\bar y$.
    Hence, $x,y$ are $\NW_F$-equivalent, so they belong to the same node.
    Since each node is forward-invariant under $F$, this shows that each node is also invariant under $F$.
\end{proof}
Next two examples show that both the compactness of the space and the openness of the map are needed for the invariance of the non-wandering set.
The first is an example of a semiflow $F$ on a two-dimensional non compact topological manifold with a non-invariant non-wandering set.
\begin{figure}
    \centering
    \includegraphics[width=0.5\linewidth]{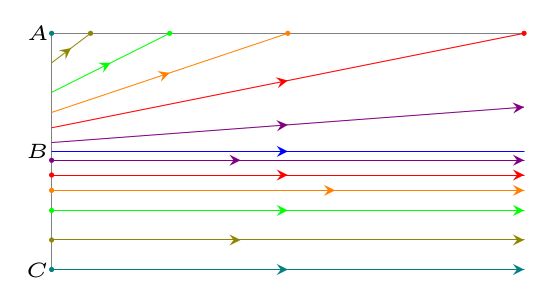}
    \caption{{\bf A semiflow with a non-invariant non-wandering set.}
    The picture shows several orbits of a semiflow $F$ on the non compact space $X$ equal to the unbounded strip shown in figure where we identify points on the horizontal half-line $h$ passing through $A$ with points on the vertical segment $BC$ so that $A$ is identified with $C$ and points going to infinity on $h$ are identified with points going to $B$ on $BC$.
    Several orbits of $F$ are shown, each one painted in a different color.
    As the picture suggests, $\cap_{t\geq0}F^t(X)=\emptyset$, i.e. no subset of $X$ is $F$-invariant.
    The non-wandering set coincides with the blue orbit.
    }
    \label{fig: nw1}
\end{figure}
\begin{example}
    Let $X$ be the quotient of the rectangle $R=[1,\infty)\times [-1,1]$ under the identification of the halfline $[1,\infty)\times\{1\}$ with the segment $\{1\}\times(0,-1]$ given by $(x,1)\sim(1,-1/x)$.
    We define on $X$ a continuous-time semiflow 
    as follows.
    For each $p\in X$, there is a point $q=(1,y)$, $y\in[0,1]$, such that $p\in\cO_F(q)$.
    For $y\in[0,1]$,
    $$
    F^t(0,y) = 
    \begin{cases}
    (1+t)(1,y),& t\leq 1/y-1;\\
    (1+t,-y),& t > 1/y.
    \end{cases}
    $$
    The action of $F$ on any other point of $X$ can be obtained from the formula above using the fact that
    $F^s(p) = F^s(F^t(q)) = F^{s+t}(q)$.
    Several orbits are illustrated in Figure~\ref{fig: nw1}, where distinct orbits are painted in different colors.
    The reader can verify that 
    $$
    NW_F=[1,\infty)\times\{0\}
    $$ 
    and that 
    $$
    F^t(NW_F)=[1+t,\infty)\times\{0\}\neq NW_F
    $$
    for every $t>0$.
    Consider now the time-1 map $f=F^1$.
    Then $NW_f=\{1,2,\dots\}\times\{0\}$ and $f(NW_f)=\{2,3,\dots\}\times\{0\}$.
    The map $f$ is open but, thanks to the fact that $X$ is not compact, not necessarily each point of $NW_f$ has a preimage. 
    In this concrete case, the point $(1,0)$ has no preimage.
\end{example}
The example above can be slightly modified to a semiflow on a compact space with non invariant non-wandering set. 
\begin{figure}
    \centering
    \includegraphics[width=0.5\linewidth]{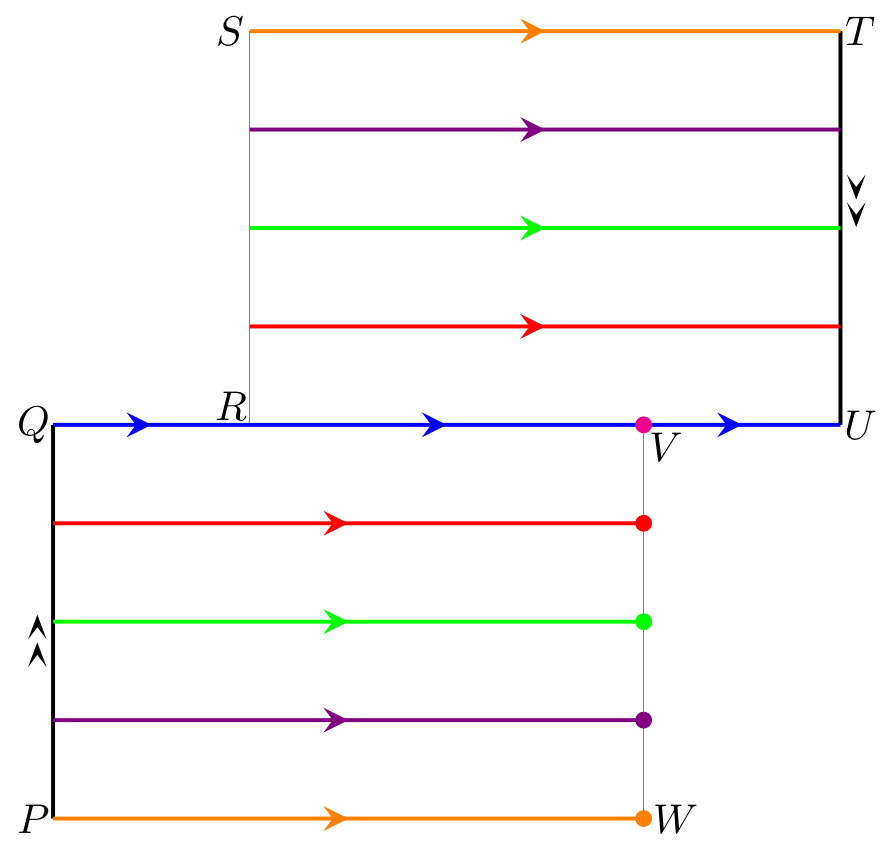}
    \caption{{\bf A semiflow with a non-invariant non-wandering set.}}
    \label{fig: nw3}
\end{figure}
\begin{example}
    Let $X$ be the Mobius strip realized as shown in Figure~\ref{fig: nw3}.
    We show in Figure~\ref{fig: nw3} several orbits of the continuous-time semiflow $F$ on $X$
    whose action on $X$ is defined as follows.
    Every orbit of $F$ is horizontal and each point, except for points on the blue one, lies on the orbit of a point on the $RS$ segment and each of these orbits asymptotes to a fixed point on the $VW$ segment. 
    The blue orbit has forward and backward limit sets equal to the magenta fixed point $V$.
    The only $F$-invariant sets of $X$ are the fixed points and the blue orbit.
    The reader can verify that each point in $RV$, including the endpoints of the segment, is non-wandering and that $NW_F=AE\cup EF$.
    Since $A$ is not fixed, $F^t(NW_F) = F^t(A)E\cup EF\neq NW_F$ for every $t>0$.
    Notice that no map $F^t$ is open at any point $p$ such that $F^t(p)=R$.
\end{example}
Assume that $F$ has compact dynamics.
Since $G_F$ contains every invariant compact set, we know from the two previous proposition that, if either $F$ is a flow or $F$ has a compact global trapping region and each $F^t$ is an open map, then $NW_F\subset G_F$.
Below we prove that actually this is always the case for a semiflow with compact dynamics.
\begin{proposition}
    \label{prop:NW_F in G}
    $NW_F\subset G_F$.
\end{proposition}
\begin{proof}
    Recall that, for every point $x$, $\Omega_F(x)\subset G_F$ and that, if $x\in NW_F$, then $\O_F(x)\cup\Omega_F(x)$ is a $\NW_F$-equivalent set.
    Hence, if $x\in NW_F\setminus G_F$, $x$ is $\NW_F$-equivalent to some point $y\in G_F$.
    It is enough, therefore, to prove that no point $x$ outside of $G_F$ can be $\NW_F$-equivalent to any point $y\in G_F$.

    Let $\eta=d(x,G_F)$ and let $\eps>0$ be small enough that $N_\eps(G_F)$ is attracted to $G_F$ and that $\eps<\eta/2$.
    Then there is a $\tau>0$ such that 
    $$ F^t(N_\eps(G_F))\subset N_\eps(G_F)\subset N_{\eta/2}(G_F)
    $$ 
    for all $t\geq\tau$.
    This means that there cannot be sequences $\eps_i\to0$, $t_i\to\infty$, $z_i\in X$ such that $d(y,z_i)<\eps_i$ and $d(x,F^{t_i}(z_i))<\eps_i$ for all $i=1,2,\dots$.
    Indeed, for almost all $i$ we have that $\eps_i<\eps$ and $t_i>\tau$. 
    Hence, for almost all $i$, we have  
    that $F^{t_i}(z_i)\in N_{\eta/2}(G_F)$, so that 
    $d(x,F^{t_i}(z_i))>\eta/2$.
    Hence, $x$ cannot be $\NW_F$-equivalent to $y$.
\end{proof}
\subsection{The Graph of \BF$\NW_F$} 
The qualitative properties of $\NW_F$ can be encoded into a graph as follows.
\begin{definition}
    The {\BF graph of $\NW_F$} (prolongational graph) is the directed graph   $\Gamma_{\NW_F}$ having the nodes of $NW_F$ as its vertices and such that there is an {\BF edge from node $N$ to node $M$} if and only if there is a $x\in N$ and a $y\in M$ such that $x\NWto y$.
    Sometimes we call the edge {\bf strong} if there is a bitrajectory $b$ such that $\alpha(b)\subset N$ and $\omega(b)\subset M$; sometimes we call an edge is {\bf weak} if it is not strong.
\end{definition}
\begin{lemma}[De Leo \& Yorke, 2025~\cite{DLY25}, Proposition 3.1.1]
    \label{prop:Down(compact)}
    Let $K\subset X$ and set 
    $
    \Up_{\NW_F}(K)=\{y:\text{ there is }x\in K\text{ such that } y\NWto x\}.$
    Then, if $K$ is compact, the set $\Up_F(K)$ is closed.
\end{lemma}
\begin{lemma}[De Leo \& Yorke, 2025~\cite{DLY25}, Proposition 4.1.5, case (3)]
    \label{lemma: case 3}
    Assume that $\Omega_F(x)\neq\emptyset$. 
    Then, if $x\NWto y$, either $y\in\cO_F(x)$ or $\Omega_F(x)\NWto y$.
\end{lemma}
\begin{proposition}
    \label{prop: F|Q}
    $\Gamma_{\NW_{F_Q}}=\Gamma_{\NW_F}$ for each fat $Q\in\cK_F$. 
\end{proposition}
\begin{proof}
    Recall that fat compact global trapping regions exist by Proposition~\ref{prop: F has a fat compact global trapping region}.
    The claim is an immediate consequence of Proposition~\ref{prop:NW_F in G} and the fact that, since $Q$ is a forward-invariant neighborhood of $NW_F$, for $\eps>0$ small enough each $(F,\eps)$-link from $x\in NW_F$ to $y\in NW_F$ lies entirely in $Q$.
\end{proof}
The example below shows that Proposition~\ref{prop: F|Q} does not extend to the restriction of $F$ to $G_F$.
%
\begin{figure}
    \centering
    \includegraphics[width=0.5\linewidth]{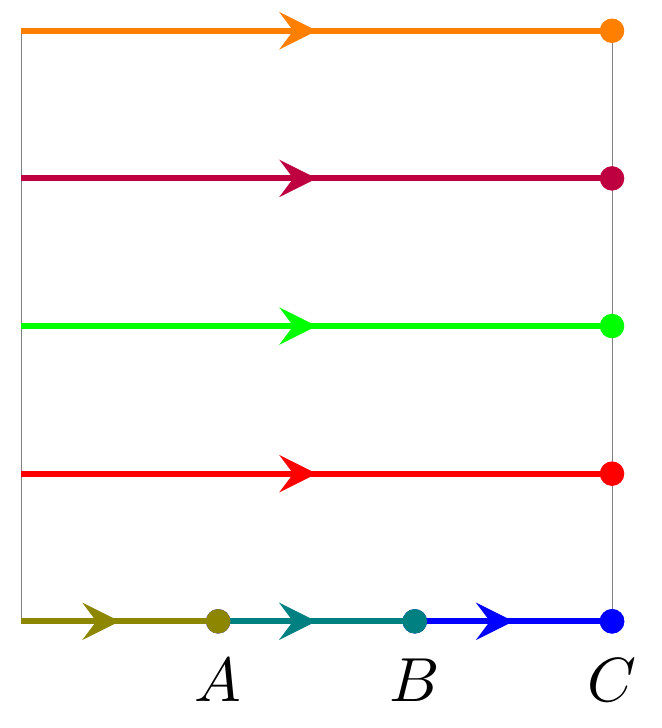}
    \caption{An example of semiflow $F$ where the prolongational graph of $F$ and that of the restriction of $F$ to its global attractor do not coincide.}
    \label{fig: prolongation graph on G_F}
\end{figure}
\begin{example}
    Consider the semiflow $F$ on the unit square sketched in Figure~\ref{fig: prolongation graph on G_F}.
    Except for the fixed points, each point moves rightwards.
    The orbit of each point $(0,y)$, $y\in(0,1]$, is the whole segment $[0,1]\times\{y\}$ and each point $(1,y)$, $y\in(0,1]$, is a fixed point.
    On the side $[0,1]\times\{0\}$ there are three fixed points $A=(1/3,0)$, $B=(2/3,0)$ and $C=(1,0)$.
    
    The reader can verify that 
    $$
    G_F = [1/3,1]\times\{0\}\bigcup\{1\}\times[0,1]
    $$
    and that no point but the fixed ones are non-wandering under $F$.
    Moreover, $\Gamma_{\NW_F}$ has three edges: one from $A$ to $B$, one from $B$ to $C$ and one from $A$ to $C$.
    The first two are strong, the last one is weak.

    Denote by $F_G$ the restriction of $F$ to $G_F$.
    Since it is $F$-invariant, $G_F$ is the global attractor of $F_{G}$ as well.
    Moreover, one can verify that $F_{G}$ and $F$ have the same non-wandering sets and the same non-wandering nodes.
    Furthermore, the edges $A\NWto B$ and $B\NWto C$ are also edges for the graph of $F_G$, since they come from two-sided trajectories that are contained in $G_F$.
    On the other side, the edge $A\NWto C$ is not in the graph of $F_{G}$ since that edge comes from orbits that are outside of $G_F$.
    Hence, the graph of $\cP_{F_{G}}$ does not coincides with the graph of $\cP_F$.
    
    
\end{example}
%
Below we give some further detail on which edges of $\Gamma_{\cP_{F}}$ might not be present in $\Gamma_{\cP_{F_{G}}}$.
\begin{proposition}
    If $\Gamma_{\cP_{F}}$ has a strong edge from node $M$ to node $N$, the same edge is present in $\Gamma_{\cP_{F_{G}}}$.
\end{proposition}
\begin{proof}
    By definition, there must be a bitrajectory $b\subset X$ such that $\alpha(b)\subset M$ and $\omega(b)\subset N$.
    Then $\alpha(b)\cup b\cup\omega(b)$ is compact and $F$-invariant and so it must be contained in $G_F$.
\end{proof}
In other words, only weak edges might not be inherited by $\Gamma_{\cP_{F_{G}}}$.
\begin{definition}
    Let $M,N$ be nodes of $\Gamma_{\NW_F}$ with $M\NWto N$, i.e. there is an edge from $M$ to $N$.
    We say $M$ and $N$ are {\BF adjacent} if, given a node $K$ with $M\NWto K$ and $K\NWto N$, either $K=M$ or $K=N$.
\end{definition}
\begin{proposition}
    \label{prop: strong edge}
    Let $M$ and $N$ be two adjacent nodes of $\Gamma_{\cP_{F}}$.
    Then the edge from $M$ to $N$ is strong.
\end{proposition}
\begin{proof}
    We need to prove that, given two nodes $M,N$ of $\cP_F$, if $M\NWto N$, then $M\NWtoFG N$.
    By Proposition~\ref{prop: F has a fat compact global trapping region}, there is a fat $Q\in\cK_F$ and, by Proposition~\ref{prop: F|Q}, $\Gamma_{\NW_{F_Q}}=\Gamma_{\NW_F}$.
    Hence, we restrict $F$ to $Q$.

    Assume first that $M\cap N\neq\emptyset$ and let $x\in M\cap N$. Then $\Omega_F(x)\subset M\cap N\subset N$.
    Since non-wandering nodes of $F$ are fully $F$-invariant, $x$ has a backward trajectory inside $M$.
    Hence, there is a bitrajectory $b$ with $\alpha(b)\subset M$ and $\omega(b)\subset N$.

    Assume now that $M$ and $N$ are disjoint. Since $Q$ is compact, $d(M,N)>0$.
    We claim that, for every $t>0$, there exist disjoint open sets $U\supset M$ and $V\supset N$ such that $F^t(U)\cap V=\emptyset$.
    If it were not so then, for every $\eps>0$, there would be a $x$ with $d(x,M)<\eps$ such that $d(F^t(x),N)<\eps$.
    By taking $\eps_n=1/n$, we can build sequences $x_n$ such that $d(x_n,M)<\eps_n$ and $d(F^{t}(x_n),N)<\eps_n$.
    We can assume without loss of generality that $x_n\to y$. 
    Clearly $y\in M$, so that $F^t(y)\in M$, but, by continuity, we should also have $d(F^{t}(y),N)=0$, which is not possible since $d(M,N)>0$.

    An important consequence of the existence of such $U$ and $V$ is that every $\eps$-$\alpha\omega$chain has at least a point that does not belong to $U\cup V$. 
    Take again $\eps_n=1/n$ and let $C_n$ be $\eps_n$-$\alpha\omega$chains from $M$ to $N$.
    Since $X$ is compact, the points of these chains have at least an accumulation point $z$ not belonging to $U\cup V$.
    Since $z$ is limit of a sequence of points whose orbit passes within $\eps_n$ from $N$, we have that, by continuity, $\Omega_F(z)\subset N$.

    Now, denote by $E$ the set of all accumulation points of all $\eps$-$\alpha\omega$chains from $M$ to $N$ for all $\eps>0$.
    Then $E$ is compact (since it is closed and a subset of $Q$) and $F$-invariant.   
    Indeed, let $x_i\to z$, with $x_i\in C_i$. Then, for each $t>0$ and $x_i\in C_i$, there is a predecessor $y_i$ in $C_i$ such that $F^t(y_i)=x_i$.
    We can assume without loss of generality that $y_i$ converge to some $w$ in $E$, so that $F^t(w)=z$.

    Hence, for every $z\in E$ we can build a backward trajectory of $z$ in $E$ and therefore a bitrajectory $b$ based at $z$.
    Let $K$ be the node such that $\alpha(b)\subset K$. By construction, each point $x\in K$ is a limit of points belonging to $\eps$-$\alpha\omega$chains from $M$ to $N$ for $\eps\to0$. 
    Hence, we can break each such chain into a $\eps$-$\alpha\omega$chain from $M$ to $K$ and another $\eps$-$\alpha\omega$chain from $K$ to $N$, so that $M\NWto K$ and $K\NWto N$. Since $M$ and $N$ are adjacent, this means that $K=M$. Hence, $b$ is a bitrajectory that runs from $M$ to $N$.   
\end{proof}
\begin{corollary}
    Any two nodes that are adjacent in $\Gamma_{\NW_F}$ are adjacent in $\Gamma_{\NW_{F_G}}$
\end{corollary}
\begin{definition}
    Let $\cN_1,\cN_2$ be two mutually disjoint collections of nodes of $NW_F$ such that each node of $NW_F$ belongs to either $\cN_1$ or $\cN_2$ and denote by $N_1$ and $N_2$ the unions of all points in, respectively, the nodes  $\cN_1$ and $\cN_2$.
    We say that $\cN_1,\cN_2$ are a {\bf nodes partition} of $\Gamma_{\NW_F}$ if $N_1\cap N_2=\emptyset$.

    We say that a prolongational graph $\Gamma_{\NW_F}$ is {\bf connected} if, for each of its nodes partitions $\cN_1,\cN_2$, there is an edge from a node of $\cN_1$ to a node of $\cN_2$ or viceversa.
\end{definition}
Notice that, given any partition $\cN_1,\cN_2$, we have that $N_1\cup N_2=NW_F$.
%
\begin{theorem}
    \label{thm:NW}
    Assume that $G_F$ is connected.
    Then $\Gamma_{\NW_F}$ is  connected.
\end{theorem}
%
%
\begin{proof}

    Let $\cN_1,\cN_2$ be a nodes partition of $NW_F$ and denote the respective sets of points by $N_1$ and $N_2$.
    Since $G_F$ is connected and $N_1\cup N_2$ is not, $E\bydef G_F\setminus(M\cup N)\neq\emptyset$.

    Let $x\in E$.
    Since $G_F$ is invariant, there is at least a bitrajectory $b$ through $x$.
    By Proposition~\ref{prop: Om(x)}, its limit sets $\alpha(x)$ and $\omega(x)$ belong to some node of $NW_F$.

    Assume first that there is a bitrajectory $b$ such that $\alpha(b)\subset N_1$ and $\omega(b)\subset N_2$.
    Then this $b$ is a strong edge from some node in $\cN_1$ to some node of $\cN_2$.

    Assume now that, for each $x\in E$, each bitrajectory through $x$ is such that $\alpha(b)\subset N_1$ and $\omega(b)\subset N_1$.
    Since $G_F$ is connected, for every $\eps>0$ there is a $x_\eps\in E\cap N_\eps(N_2)$.
    Then, since by the working hypothesis $\Omega_F(x_\eps)\subset N_1$, the trajectory of $x_\eps$ can be used to build a $(F,\eps)$-link, i.e. a weak edge, from some node of $\cN_2$ to some node of $\cN_1$.

    Assume finally that, for each $x\in E$, each bitrajectory through $x$ is such that either $\alpha(b)\subset N_1$ and $\omega(b)\subset N_1$ or $\alpha(b)\subset N_2$ and $\omega(b)\subset N_2$.
    Let $E_1$ be the set of all points of the first type and $E_2$ the set of points of the second type.
    Then $N_1\cup E_1$ and $N_2\cup E_2$ are closed and disjoint and their union is $G_F$, contradicting the hypothesis that $G_F$ is connected. 

    Every other case reduces to one of the three cases above. 
    Hence, in any case, for every partition $\cN_1,\cN_2$ of $NW_F$ we have an edge between a node of $\cN_1$ and a node of $\cN_2$, namely $\Gamma_{\NW_F}$ is connected.
    %
\end{proof}
\section{Streams}
\label{sec: streams}
\begin{definition}
    We call {\BF $F$-stream} (or simply {\bf stream}, when there is no ambiguity) on $X$ a closed quasi-order $S$ that is an extension of $\cO_F$.
    When $(x,y)\in S$, we use the notation {\BF $x\Sto y$}.
    We write {\BF$x\Seq y$}, and we say that $x$ and $y$ are {\BF$S$-equivalent}, if $x\Sto y$ and $y\Sto x$.
    We say that $x$ is {\BF $S$-recurrent} if either $x$ is fixed or there is a $y\neq x$ such that $x\Seq y$.
    We denote by $\cR_S$ the set of all $S$-recurrent points.
    We say that a set $M\subset\sR_{S}$ is $S$-equivalent if all points of $M$ are mutually $S$-equivalent.
    We call {\bf nodes} the equivalence classes of $\cR_{S}$ with respect to $\Seq$.
    If $x\in\cR_S$, we denote by $\Node_S(x)$ the node containing $x$.
\end{definition}
Recall that, as implicitly assumed in the definition above, $\Seq$ is an equivalence relation.
\begin{proposition}
    Any intersection of $F$-streams is a $F$-stream.
\end{proposition}
Hence, each semiflow $F$ has a largest stream, the relation $X\times X$, and a smallest stream, the intersection of all of its streams, all of which contain the prolongational relation $\cP_F$.
See more about the smallest stream in Section~\ref{sec: A_F}.
\begin{proposition}
    \label{prop:streams}
    Let $S$ be a $F$-stream. Then:
    \begin{enumerate}
        \item For each $x$, $\Omega_F(x)$ is a $S$-equivalent set. In particular, there is a node $N$ of $S$ such that $\Omega_F(x)\subset N$.
        \item If $x$ is $F$-recurrent, then $\cO_F(x)\subset\Omega_F(x)\subset\Node_S(x)$.
    \end{enumerate}
\end{proposition}
The qualitative properties of a stream can be encoded in a graph as follows.
\begin{definition}
    \label{def: N1>=N2}
    Let $S$ be a $F$-stream.
    Given two sets $A,B\subset X$, we write $A\Sto B$ if and only if $x\Sto y$ for each $x\in A$ and $y\in B$.
    Given a $F$-stream, the graph of $S$, denoted by {\bf\BF$\Gamma_S$}, is the directed graph whose nodes are the nodes of $\sR_S$ and such that there is an edge from a node $N_1$ to a node $N_2\neq N_1$ if and only if $N_1\leadsto N_2$.
    A graph $\Gamma_S$ is {\bf connected} if, whenever $S=C_1\cup C_2$, with $C_1$ and $C_2$ closed and disjoint sets each of which is union of nodes of $S$, there is an edge from a node of $C_1$ to a node of $C_2$ or viceversa.
\end{definition}
\begin{theorem}
    \label{thm: connectedness}
    Assume that $G_F$ is connected and 
    let $S$ be a stream of $F$. 
    Then $\Gamma_S$ is connected.  
\end{theorem}
\begin{proof}
    Each $F$-stream $S$ has two types of nodes: those that are an extension of non-wandering nodes of $F$ and those that are not. 
    In this proof, we will refer to the first type as ``$\Omega$-nodes''.
    Notice that $S$-nodes that are not $\Omega$-nodes cannot be forward-invariant under $S$.
    If they were, indeed, they would contain a $F$-recurrent point, since every node is compact under the theorem's hypotheses, and every $F$-recurrent is in some non-wandering node.
    Ultimately, each $S$-node either is an $\Omega$-node or has an edge from itself to at least one $\Omega$-node.

    Suppose now that $S=C_1\cup C_2$ with $C_1,C_2$ closed and disjoint sets each of which is a union of nodes.
    If either one of the two, say $C_1$, only contains nodes that are not $\Omega$-nodes, then there is at least an edge from $C_1$ to $C_2$ because the limit sets of all points of $C_1$ lie in $C_2$.
    Suppose now that both $C_1$ and $C_2$ contain $\Omega$-nodes. 
    Then, since $\Gamma_{\NW_F}$ is connected by Theorem~\ref{thm:NW}, there is at least an edge between an $\Omega$-node in $C_1$ and an $\Omega$-node in $C_2$.
    Otherwise, it would be possible to sort the non-wandering nodes into two disjoint closed sets so that there would be no edge between the two sets, namely $\Gamma_{\NW_F}$ would not be connected.
    Hence, $\Gamma_S$ is connected.
\end{proof}

\medskip\noindent
{\BF $\Omega$streams.}
An immediate consequence of stream's transitivity is that, for all $y\in\cO_F(x)$,
\beqn
\label{eq: Ostreams}
\Down_{\NW_F}(x) \supset \cO_F(x)\cup\Down_{\NW_F}(y).
\eeqn
It turns out that streams for which the sets at the left and right hand sides above are equal enjoy rather special properties.
In this section we illustrate some of them.

\begin{definition}
    We say that a $F$-stream $S$ is a {\BF $\Omega$stream} if
    \beq
    \Down_S(x) = \cO_F(x)\cup\Down_S(y)
    \eeq
    for every $y\in\cO_F(x)$.
\end{definition}
\begin{proposition}
    A $F$-stream $S$ is a $F$-$\Omega$stream if and only if    
    \beq
    \Down_S(x) = \cO_F(x)\cup\Down_S(\Om_F(x)).
    \eeq
\end{proposition}
\begin{definition}
    We say that a node $N$ of an $F$-stream $S$ is {\em dynamical} if $N$ contains $F$-recurrent points. 
\end{definition}  
\begin{proposition}[De Leo \& Yorke, 2025~\cite{DLY25}, Proposition 5.3.7]
    \label{prop:minimal streams}
    Let $S$ be an $F$-$\Omega$stream. Then:
    \begin{enumerate}
        \item If $x$ is $S$-recurrent, then $\cO_F(x)\subset\Node_S(x)$.
        \item $\Down_S(\Om_F(x))\cap\sR_S=\Down_S(x)\cap\sR_S$.
        \item Every node of $S$ is closed and forward-invariant under $F$.
        \item Every node of $S$ is dynamical.
        \item $\cR_S$ is closed and forward-invariant under $F$.
        \item $S$ is forward-invariant under the natural action induced by $F$ on $X\times X$.
        \item $\Down_S(M)$ is forward-invariant under $F$ for each $M\subset X$.
        \item $\Up_S(C)$ is forward-invariant under $F$ for each set $C\subset\cR_S$ that is union of nodes of $S$.
    \end{enumerate}
\end{proposition}
%
%
%
\subsection{Trapping regions of streams} 
%
Recall that a trapping region for $F$ is a closed set $Q$ such that $F^t(Q)\subset Q$ for all $t\geq0$.
By analogy, we provide the following definition in case of streams.
\begin{definition}
    Given a stream $S$, we say that a closed set $Q$ is a trapping region for $S$ 
    if $\Down_S(Q)\subset Q$.
\end{definition}
Notice that, since $\cO_F\subset S$, each trapping region for a $F$-stream $S$ is also a trapping region for $F$.

%
\begin{lemma}[De Leo \& Yorke, 2025~\cite{DLY25}, Lemma 5.5.2]
    \label{lemma:S trapping region}
    Let $S$ be a $F$-stream on $X$ and let $Q$ be a trapping region for $S$.
    The following hold:
    \begin{enumerate}
        \item if $S'$ is a substream of $S$, then $Q$ is a trapping region for $S'$;
        \item if a node $N$ of $S$ has some point in $Q$, then $N\subset Q$;
        \item there is no edge in $\Gamma_S$ from any node in $Q$ to any node outside $Q$;
        \item if $S$ is a $\Omega$stream, then $x\Sto y$, with $x,y\not\in Q$, if and only if $x\Fto y$;
        \item if $S$ is a $\Omega$stream, $\cR_S\subset Q$.
    \end{enumerate}
\end{lemma}
%
%
\begin{proposition}
    Let $\cS_S$ be the set of substreams of a $\Omega$stream $S$ and let $Q$ be a trapping region for $S$.
    Then $S_1,S_2\in\cS_S$ coincide if and only if their restriction to $Q$ coincide.
\end{proposition}
\begin{proof}
    By hypothesis, $S_1\cap(Q\times Q)=S_2\cap(Q\times Q)$.
    We need to prove that $(x,y)\in S_1$ if and only if $(x,y)\in S_2$ when either $x\not\in Q$ or $y\not \in Q$ (or both).
    
    Suppose that $(x,y)\in S_1$.
    By the proposition above (point 4), when both $x,y$ are outside of $Q$, $(x,y)$ belongs to all substreams of $S$.
    If $x\in Q$, since $Q$ is a trapping region for $S$, then we must have that $y\in Q$.
    The last case is when $x\not\in Q$ and $y\in Q$. 
    In this case, since $S$ is a $\Omega$stream, either $(x,y)\in\cO_F$, in which case it belongs to all streams, or $y\in\Down_{S_1}(\Om_F(x))$.
    Since $S_1$ and $S_2$ coincide inside $Q$, then $(x,y)\in S_2$.
\end{proof}
\subsection{Streams of chains}
The most important generalizations of recurrent points in literature are built out of chains, as defined below.
These are Auslander's generalized recurrent points~\cite{Aus63}, Easton's strong chain-recurrent points~\cite{Eas77} and Conley's chain-recurrent points~\cite{Con72}.
In this section we define corresponding streams, which we call the Auslander stream, the $\Sigma$chains stream and the chains stream respectively, that have, respectively, the sets above as their set of recurrent points.
\subsubsection{Discrete-time chains streams}
\label{sec:CR}
Recall that, given a map $f$ on $X$, a $(f,\eps)$-chain from $x$ to $y$ of length $n+1$ is a sequence $(c_0,\dots,c_n)$ such that $c_0=x$, $c_n=y$ and $d(f(c_i),c_{i+1})<\eps$ for all $i=0,\dots,n-1$.

The following technical lemma will be used several times in the rest of the section.
\begin{lemma}[Hurley, 1991~\cite{Hur91}, Lemma 1.2]
    \label{lemma: UC discrete}
    Let $f:X\to X$ be continuous  and let $x\in X$.
    Then, for any $\eps>0$ and any integer $p>0$, there exists a $\delta>0$ such that:
    \begin{enumerate}
        \item 
        $
        d(x,y)<\delta \implies d(f^k(x),f^k(y))<\eps\;\text{ for every }\;k=0,\dots,p;
        $
        \item 
        for every $(f,d,\delta)$-chain $C=(c_0,c_1,\dots,c_n)$ based at $x$ with $n\leq p$, the chain $C'$ consisting in the pair of points $c_0,c_n$ is a $(f^n,d,\eps)$-chain.
        Equivalently, 
        $$
        d(f^n(c_0),c_n)<\eps.
        $$
        If $f$ is uniformly continuous, then $\delta$ can be chosen independently on $x$.
    \end{enumerate}
\end{lemma}
\begin{proof}
    (1)
    We leave it to the reader.

    (2)
    The case $p=1$ is a tautology and holds for every $\delta<\eps$. 
    Assume now that the claim holds for all $p\leq p_0$.
    Since we already know that the claim is true for all $n=1,\dots,p_0$, it is enough to prove that, for every $\eps>0$, there is a $\delta>0$ such that $d(f^{p_0+1}(c_0),c_{p_0+1})<\eps$
    
    Consider a chain $(c_0,\dots,c_{p_0+1})$ and notice that 
    $$
    d(f^{p_0+1}(c_0),c_{p_0+1})
    \leq
    d(f^{p_0+1}(c_0),f(c_{p_0}))
    +
    d(f(c_{p_0}),c_{p_0+1}).
    $$
    By continuity, there is an $\eta>0$ such that $d(y,z)<\eta$ implies $d(f(y),f(z))<\eps/2$.
    By the inductive assumption, 
    there is a $\delta_1>0$ such that, if $(c_0,\dots,c_{p_0})$ is a $(f,d,\delta_1)$-chain, then $d(f^{p_0}(c_0),c_{p_0})<\eta$.
    Therefore, if $(c_0,\dots,c_{p_0})$ is a $(f,d,\delta_1)$-chain, then $d(f^{p_0+1}(c_0),f(c_{p_0}))<\eps/2$.
    
    Assume now that $(c_0,\dots,c_{p_0+1})$ is a $\eps/2$-chain.
    Then $d(f(c_{p_0}),c_{p_0+1})<\eps/2$.
    Finally, let $\delta=\min\{\eps/2,\delta_1\}$.
    Then, for every $\delta$-chain based at $x$, $$d(f^{p_0+1}(c_0),c_{p_0+1})<\eps/2+\eps/2=\eps.$$
\end{proof}
\begin{corollary}[Hurley, 1991]
    \label{cor: chains with infinite length}
    Let $y\in\Down_{\cC_{F,d}}(x)$ and $y\not\in\cO_F(x)$.
    Then, if $\eps_n\to0^+$ and $C_n$ is a $(F,d,\eps_n)$-chain from $x$ to $y$, the length of the $C_n$ diverges as $\eps_n$ goes to 0.
\end{corollary}

\begin{definition}
    Given a discrete-time semi-flow $F$ on $X$ and a metric $d$ compatible with the topology of $X$, we call {\BF $(F,d,\eps)$-chains stream} the relation
    \beq
    \cC_{F,d,\eps} = \{(x,y):\text{ there is a $(F,d,\eps)$-chain from $x$ to $y$}\}.
    \eeq
    We call {\BF $(F,d)$-infinitesimal chains stream} (or simply {\bf chains stream}) the relation
    \beq
    \cC_{F,d} = \bigcap_{\eps>0}\cC_{F,d,\eps}.
    \eeq
\end{definition}
%
%
\begin{proposition}[De Leo \& Yorke, 2025~\cite{DLY25}, Proposition 6.1.2]
    \label{prop: C is a Ostream}
    $\cC_{F,d,\eps}$ is a stream for every $\eps>0$.
    $\cC_{F,d}$ is a $\Omega$stream.
\end{proposition}
It is well known that, as we show below, when $X$ is compact, $\cC_{F,d}$ does not depend on the metric $d$.
\begin{proposition}
    \label{lemma: C1=C2}
    Let $X$ be compact and let $d_1,d_2$ be any two metrics generating the topology of $X$.
    Then $\cC_{F,d_1}=\cC_{F,d_2}$ for any discrete-time semiflow $F$ on $X$.
\end{proposition}
\begin{proof}
    Suppose that there are $x,y\in X$ such that $y$ is $\cC_{F,d_1}$-downstream but not $\cC_{F,d_2}$-downstream from $x$.

    Let $C_i$ be a sequence of $(d_1,\eps_i)$-chains from $x$ to $y$, with $\eps_i\to0$. 
    Since $y$ is not $\cC_{F,d_2}$-downstream from $x$, there is a $\delta>0$ such that, for each $i$, there is a point $x_{k_i}$ on $C_i$ such that
    $
    d_2(f(x_{k_i-1}),x_{k_i}) > \delta.
    $

    Since $X$ is compact, we can assume without loss of generality that these $x_i$ converge to a point $z$. 
    Since $d_1$ and $d_2$ are equivalent, $x_{k_i}\to z$ for both $d_1$ and $d_2$. This means that, for every $\eta>0$, we can find an $i$ such that
    $$
    d_1(f(x_{k_i-1}),x_{k_i}) < \eta,\;
    d_1(z,x_{k_i}) < \eta,\;
    d_1(z,f(x_{k_i-1})) < \eta.
    $$
    In particular, also $f(x_{k_i-1}) \to z$ with respect to the $d_1$ distance, and so it does with respect to $d_2$ as well.
    This means that, for every $\eta > 0$, we can find an $i$ large enough such that:
    $$
    d_2(z,x_{k_i}) < \eta,
    d_2(z,f(x_{k_i-1})) < \eta.
    $$
    On the other side, we also have (see above) that
    $d_2(f(x_{k_i-1}),x_{k_i}) > \delta$.
    These three inequalities are incompatible with the triangular inequality for $\eta$ small enough.
    Hence, we must have that $\cC_{F,d_1}=\cC_{F,d_2}$. 
\end{proof}
Next example, that extends an example by Alongi and Nelson in~\cite{AN07}, shows that this is not the case when $X$ is not compact.
\begin{example}
    Let $X=\{(x,y):\;x\in\R,\;y\geq1\}$ and let $F$ be the discrete-time flow of the map
    $$
    f(x,y)=(x+1,y).
    $$
    Recall that the upper half-plane with the Riemannian metric tensor $(dx^2+dy^2)/y^2$ is a model of hyperbolic geometry; we denote the corresponding distance function by $d_H$ and recall that
    $$
    d_H((x,y_1),(x,y_2)) = |y_1-y_2|,\;
    d_H((x_1,y),(x_2,y)) = \frac{|x_1-x_2|}{y}.
    $$
    Denote finally by $d_E$ the euclidean distance function.
    It is easy to verify that $\cC_{F,d_E}=\cO_F$, the smallest possible stream on $X$, and that $\cR_{\cC_{F,d_E}}=\emptyset$.
    On the other side, we show below that $\cC_{F,d_H}=X\times X$, the largest possible stream on $X$.

    Let $\eps>0$.
    The strategy to build a $(F,\eps)$-chain from
    any point of $X$ to any other one
    is the following. 
    First, it is possible to make the $y$ coordinate of the chain's elements arbitrarily large by adding $\eps/2$ at every step:
    $$
    d_H(f(x,y),(x+1,y+\eps/2))=\eps/2<\eps.
    $$
    Once $y$ is larger than $2/\eps$, it is possible to make the $x$ coordinate move by an arbitrary amount against the flow by subtracting at most 1 to it at every step:
    $$
    d_H(f(x,y),(x-1,y))<\frac{2}{2/\eps}<\eps.
    $$
    If needed, the point can be then lowered by repeating the first step but now subtracting $\eps/2$ at every step.
    Since the flow naturally moves horizontally points rightwards, it is clear that every point of $X$ is $\cC_{F,d_H}$-downstream from any other point of $X$.
    In particular, $\cR_{\cC_{F,d_H}}=X$ and the graph of $\cC_{F,d_H}$ has a single node and no edge.
\end{example}
%
Below we show that, when $F$ has compact dynamics, $\cC_{F,d}$ is purely topological.
\begin{proposition}
    \label{prop: G_F is a trapping region}
    $G_F$ is a trapping region for $\cC_{F,d}$.
\end{proposition}
\begin{proof}
    We need to prove that $\Down_{\cC_{F,d}}(G_F)\subset G_F$.
    Let $(x,y)\in\cC_{F,d}$ with $x\in G_F$ and $y\not\in G_F$ and set $\rho=d(y,G_F)$.
    Notice that $\rho>0$ since $G_F$ is compact.
    We can assume without loss of generality that $N_\rho(G_F)$ is compact. 
    
    By Proposition~\ref{prop: F has a fat compact global trapping region}, there is a fat $Q\in\cK_F$ such that $Q\subset N_{\rho/2}(G_F)$.
    Let $\eps\in(0,\rho)$ such that $N_\eps(G_F)\subset Q$.
    There is an integer $N>0$ such that $F^t(N_\eps(G_F))\subset N_{\eps/2}(G_F)$ for $t\geq N$.
    By Lemma~\ref{lemma: UC discrete}, there is $\eta\in(0,\eps)$ such that every $(F,d,\eta)$-chain $(c_0,\dots,c_n)$, $n\leq N$, with $c_0\in N_\rho(G_F)$ satisfies $d(F^n(c_0),c_n)<\eps/2$. 

    Let now 
    $
    (c_0,\dots,c_N,\dots,c_{2N},\dots,c_{kN},\dots,c_{kN+d}), d<N,
    $
    be a $(F,d,\zeta)$-chain from $x$ to $y$ (i.e. $c_0=x$ and $c_{kN+d}=y$) with $\zeta\leq\eta$.
    Since $c_0=x\in G_F$ and $G_F$ is invariant, $d(F^N(c_0),c_N)<\eps/2$, namely $c_N\in N_{\eps/2}(G_F)$.
    Since $c_N\in N_{\eps}(G_F)$, then $F^N(c_N)\in N_{\eps/2}(G_F)$ and $d(F^N(c_N),c_{2N})<\eps/2$, so that $c_{2N}\in N_{\eps}(G_F)$.
    By repeating this argument a finite number of times we find that $c_{kN}\in N_{\eps}(G_F)$ and that $d(F^d(c_{kN}),y)\leq\eps/2$.
    Since $N_{\eps}(G_F)\subset Q$ and $Q$ is forward-invariant, then $F^d(c_{kN})\in Q\subset N_{\rho/2}(G_F)$ and so 
    $$
    d(y,G_F)\leq d(y,F^d(c_{kN})) + d(F^d(c_{kN}),G_F) \leq \eps/2+\rho/2<\rho,
    $$
    contradicting the initial hypothesis that $d(y,G_F)=\rho$.
    Hence, no point outside $G_F$ can be $\cC_{F,d}$-downstream from a point of $G_F$.
\end{proof}

\begin{proposition}
    \label{prop: d1=d2}
    Let $F$ be a semi-flow with compact dynamics and let $d_1,d_2$ be any two equivalent metrics on $X$.
    Then $\cC_{F,d_1}=\cC_{F,d_2}$.
\end{proposition}
\begin{proof}
    By Proposition~\ref{prop: G_F is a trapping region}, $G_F$ is a trapping region for both $\cC_{F,d_1}$ and $\cC_{F,d_2}$. 
    By Proposition~\ref{prop: C is a Ostream}, both $\cC_{F,d_1}$ and $\cC_{F,d_2}$ are $\Omega$streams.
    We now show that they are identical.
    
    Consider first the case $x,y\not\in G_F$.
    By case 4 of Lemma~\ref{lemma:S trapping region}, $(x,y)\in\cC_{F,d_i}$, $i=1,2$, if and only if $y\in\cO_F(x)$, which is a condition independent on the distance used on $X$.
    Hence, $\cC_{F,d_1}$ and $\cC_{F,d_2}$ agree on pairs of points outside of $G_F$.

    The case $(x,y)\in\cC_{F,d_i}$ with $x\in G_F$ and $y\not\in G_F$ cannot happen for either $i=1$ or $i=2$ because of Proposition~\ref{prop: G_F is a trapping region}.
    
    Consider now the case where $x\not\in G_F$ and $y\in G_F$.
    Since $\cC_{F,d_1}$ and $\cC_{F,d_2}$ are $\Omega$streams, either $y\in\Omega_F(x)$ or there are $z_1,z_2\in\Omega_F(x)$ such that $x\CRtoI z_1\CRtoI y$ and $x\CRtoII z_1\CRtoII y$.
    Since each pair $(x,y)$ with $y\in\Omega_F(x)$ is in both $\cC_{F,d_1}$ and $\cC_{F,d_2}$,
    we are now left with the case when $x,y\in G_F$.
    The fact that $(x,y)\in\cC_{F,d_1}$ if and only if $(x,y)\in\cC_{F,d_2}$ is proven in Proposition~\ref{lemma: C1=C2}.
\end{proof}
Because of the proposition above, from now on we will denote the chains stream of $F$ by simply $\cC_F$.
\begin{proposition}[Norton, 1995~\cite{Nor95}]
    \label{prop: Norton}
    Let $X$ be compact and let $F$ be a discrete semiflow on $X$.
    Then $\cR_{\cC_F}$ is $F$-invariant.
\end{proposition}
The following proposition generalizes to our setting Douglas Norton's result above. 
\begin{proposition}
    \label{prop: F-inv}
    $\cR_{\cC_{F}}$ and all of its nodes are compact and $F$-invariant and $\cR_{\cC_{F}}\subset G_F$.
\end{proposition}
\begin{proof}
    This is an immediate consequence of Proposition~\ref{prop: G_F is a trapping region} and Proposition~\ref{prop: Norton}.
\end{proof}
\begin{lemma}
    \label{lemma: K}
    Let $M,N$ be distinct nodes of $\cC_F$ and assume that $M\CRto N$.
    Then there exists a compact set $K\subset X$ such that:
    \begin{enumerate}
        \item $K$ is $F$-invariant;
        \item $M\CRKto N$.
    \end{enumerate}
\end{lemma}
\begin{proof}
    Set $f=F^1$ and let $x\in M$ and $y\in N$.
    By hypothesis, for every $\eps>0$, there is a $(F,\eps)$-chain from $x$ to $y$.
    Let $\eps_i\to0^+$ and let $C_i=(c_{i,0},\dots,c_{i,n_i})$ be a $(F,\eps_i)$-chain in $X$ from $x$ to $y$.
    Each $C_i$ is a finite sequence of points and so is compact.
    The argument used to prove Proposition~\ref{prop: G_F is a trapping region} shows that, given any $\eta>0$, for $\eps_i$ small enough, $C_i\subset N_\eta(G_F)$.
    Hence, we can assume without loss of generality that each $C_i$ lies in some compact neighborhood $U$ of $G_F$.
    Recall that the set of all compact subsets of a compact space is complete (with respect to the Hausdorff metric).
    Hence, the $C_i$, possibly after passing to a subsequence and relabeling the indices, converge to some compact set $C\subset U$.

    We claim that $K=C\cup M$ satisfies the properties in the statement.
    
    {(1) $K$ is $F$-invariant.}
    
    By Proposition~\ref{prop: F-inv}, $M$ is $F$-invariant and so it is enough to consider the case when $z\in C\setminus M$.
    By construction, there is a sequence $x_i\in C_i$ such that $x_i\to z$.
    Since $z\not\in M$, almost all $x_i$ are not the first element of $C_i$ and so they have a predecessor $y_i$ such that $d(f(y_i),x_i)<\eps_i$.
    Possibly passing to a subsequence, $y_i\to w\in K$.
    Then, by continuity of $f$ and $d$, we have that $d(f(w),y)\leq0$, namely $f(w)=y$.
    
    (2) $M\CRKto N$.

    Let $\eps>0$.
    Since $f$ is uniformly continuous in $U$, there is a $\delta>0$ such that $d(f(z),f(w))<\eps$ for every $w,z\in U$ with $d(w,z)<\delta$.
    Let $\bar n$ be such that $\eps_{\bar n}<\eps$ and $d(C_{\bar n},K)<\delta$.
    This last condition entails that, for each element $c_{\bar n,i}$ of $C_i$, there is a $z\in K$ with $d(c_{\bar n,i},z)<\delta$.
    We can assume without loss of generality that $\delta<\eps$.
    
    We set $c_0=c_{\bar n,0}=x\in K$.
    By construction, $d(f(c_{\bar n,0}),c_{\bar n,1})<\delta<\eps$.
    We set $c_1$ to the point of $K$ closest to $c_{\bar n,1}$.
    Then
    $$
    d(f(c_0),c_1) \leq d(f(c_0),c_{\bar n,1}) + d(c_{\bar n,1},c_1)<2\eps.
    $$
    Now, we set $c_2$ to the point of $K$ closest to $c_{\bar n,2}$.
    Then
    $$
    d(f(c_1),c_2) 
    \leq 
    d(f(c_1),f(c_{\bar n,1})) 
    + 
    d(f(c_{\bar n,1}),c_{\bar n,2})
    +
    d(c_{\bar n,2},c_2)
    <3\eps.
    $$
    By repeating this construction for each element of $C_{\bar n}$, we end up building a $(F,3\eps)$-chain in $K$ from $x$ to $y$.
    Since this can be done for every $\eps>0$, then $M\CRKto N$.
\end{proof}
Next theorem shows that all qualitative chain-recurrent properties of a semiflow with compact dynamics on a locally compact space are contained in its global attractor.
\begin{proposition}[C.~Conley, 1977~\cite{Con78}; see also C.~Robinson~\cite{CRob77}]
    \label{prop: Conley}
    Let $X$ be compact and denote by $R$ the restriction of $F$ to $\cR_{\cC_F}$. 
    Then $\cR_{\cC_{R}}=\cR_{\cC_F}$.
\end{proposition}
Notice that Conley claimed the proposition above for flows on compact spaces but his proof, as well as the one provided by Robinson and Franke in~\cite{CRob77}, works without changes for semiflows with compact dynamics.
%
%
\begin{theorem}
    \label{thm: CR G}
    Let $F_G$ be the restriction of $F$ to its global attractor.
    Then:
    \begin{enumerate}
        \item $\cR_{\cC_{F_G}}=\cR_{\cC_{F}}$;
        \item $N\subset\cR_{\cC_{F}}$ is a node of $\cC_{F}$ if and only if it is a node of $\cC_{F_G}$;
        \item $\Gamma_{\cC_{F_G}}=\Gamma_{\cC_F}$.
    \end{enumerate}
\end{theorem}
\begin{proof}
    (1)
    Since $G_F$ is a compact invariant set and $\cR_{\cC_F}\subset G_F\subset X$, then Proposition~\ref{prop: Conley} implies that both $\cR_{\cC_{F_G}}$ and $\cR_{\cC_{F}}$ are equal to the chain-recurrent set of the restriction of $F$ to $\cR_{\cC_{F}}$.
    
    (2)
    Denote by $R$ the restriction of $F$ to $\cR_{\cC_{F}}$.
    The argument used by Robinson and Franke, whose pattern we use to prove Lemma~\ref{lemma: K}, implies that, if two points are $\cC_F$-equivalent, then they are also $\cC_R$-equivalent, and so also $\cC_{F_G}$-equivalent.
    Hence, two points are $\cC_F$-equivalent if and nly if they are $\cC_{F_G}$-equivalent.

    (3)
    We know from (2) that $F$ and $F_G$ have the very same nodes.
    Now, assume that $M\CRto N$ and let $x\in M$ and $y\in N$.
    By Lemma~\ref{lemma: K}, there is a $F$-invariant compact set $K\subset X$ such that $x\CRKto y$.
    By Proposition~\ref{prop: char of G}, $K\subset G_F$ and so $x\CRGto y$.
    Hence, there is an edge from $M$ to $N$ in $\Gamma_{\cC_{F_G}}$ if and only if there is one in $\Gamma_{\cC_F}$.
\end{proof}
%
%
\subsubsection{Continuous-time chains streams}
%
So far, we only considered the case of discrete-time chains.
Here, we prove that this can be done without loss of generality because the time-1 map $f=F^1$ of a continuous-time semi-flow $F^t$ completely determines the nodes and edges of the graph of $F^t$.
Our results extend, within a compact dynamics context, the following important result by Mike Hurley:
%
\begin{thmE}[Hurley, 1995~\cite{Hur95}]
    Let $F$ be a continuous-time semi-flow on a compact metric space $X$ and let $f=F^1$ be the corresponding time-1 discrete-time semi-flow.    
    Then $\cR_{\cC_{F}}=\cR_{\cC_{f}}$.
\end{thmE}
Notice that the result above is not stated explicitly in~\cite{Hur95} but is rather a corollary of a more general result (Thm.~5 in~~\cite{Hur95}) that holds, in general metric spaces, 
for a stronger version of chain-recurrence, where the ``$\eps$'' of an $\eps$-chain is not a constant but rather a strictly positive function. 
In case of a compact metric space, this general result reduces to Theorem~E.

We start with the following definitions.
\begin{definition}
    Given a continuous-time semi-flow $F$ on $X$ and a metric $d$ compatible with the topology of $X$, given $\eps>0$ and $T>0$, a {\BF$(F,d,\eps,T)$-chain} of length $n+1$ from $x$ to $y$ is a sequence of $n+1$ points $c_0,\dots,c_n$ together with a finite sequence of positive real numbers $t_0,\dots,t_{n-1}$ such that:
    \begin{enumerate}
        \item $c_0=x$, $c_n=y$;
        \item $d(F^{t_i}(c_i),c_{i+1})\leq\eps$ for all $i=0,\dots,n-1$;
        \item $t_i\geq T$ for all $i=0,\dots,n-1$.
    \end{enumerate}
\end{definition}
The following technical lemma is a continuous-time analogue of Lemma~\ref{lemma: UC discrete}.
\begin{lemma}[Hurley, 1995]
    \label{lemma: UC cont}
    Let $F$ be a continuous-time semi-flow on $X$.
    Then, for any $\eps>0$, $T>0$ and $p>0$, there exists a $\delta>0$ such that:
    \begin{enumerate}
        \item 
        $
        d(x,y)<\delta \implies d(F^t(x),F^t(y))<\eps\;\text{ for every }\;t\in[0,T];
        $
        \item 
        for every $(F,d,\delta,T)$-chain $C$ with $p+1$ points $c_0,\dots,c_p$ and times $t_0,\dots,t_{p-1}$, the chain $C'$ with points $c_0,c_p$ and time $\tau=\sum_{i=0}^{p-1}t_i$ is a $(F,d,\eps,\tau)$-chain.
        Equivalently, 
        $$d(F^{\tau}(c_0),c_p)<\eps.$$
    \end{enumerate}
    If $f$ is uniformly continuous, then $\delta$ can be chosen independently on $x$.
\end{lemma}
\begin{definition}
    We call {\BF $(F,d,\eps,T)$-chains stream} the relation
    \beq
    \cC_{F,d,\eps,T} = \cO_F\cup\{(x,y):\text{ there is a $(F,d,\eps,T)$-chain from $x$ to $y$}\}.
    \eeq
    We call {\BF $(F,d,T)$-infinitesimal chains stream}
    the relation
    \beq
    \cC_{F,d,T} = \bigcap_{\eps>0}\cC_{F,d,\eps,T}.
    \eeq
    We call {\BF $(F,d)$-infinitesimal chains stream} (or simply {\bf chains stream}) the relation
    \beq
    \cC_{F,d} = \bigcap_{T>0}\cC_{F,d,T}.
    \eeq
    By analogy, given a discrete-time semi-flow $f$, we call {\BF $(f,d,\eps,N)$-chain} of length $n+1$ from $x$ to $y$ a sequence of $n+1$ points $c_0,\dots,c_n$ together with a finite sequence of positive integers $k_0,\dots,k_{k-1}$ such that:
    \begin{enumerate}
        \item $c_0=x$, $c_n=y$;
        \item $d(f^{k_i}(c_i),c_{i+1})\leq\eps$ for all $i=0,\dots,n-1$;
        \item $k_i\geq N$ for all $i=0,\dots,n-1$.
    \end{enumerate}
    We call {\BF $(f,d,\eps,N)$-chains stream} the relation
    \beq
    \cC_{f,d,\eps,N} = \cO_f\cup\{(x,y):\text{ there is a $(f,d,\eps,N)$-chain from $x$ to $y$}\}
    \eeq
    and {\BF $(f,d,N)$-infinitesimal chains stream}
    the relation
    \beq
    \cC_{f,d,N} = \bigcap_{\eps>0}\cC_{f,d,\eps,N}.
    \eeq
\end{definition}
%
The reader can verify that all the relations above are indeed streams.
\begin{proposition}
    \label{prop: cont Ostream}
    Let $F$ be a continuous-time semiflow.
    Then $\cC_{F,d}$ is an $\Omega$stream and, for every $T>0$, $\cC_{F,d,T}$ is an $\Omega$stream.
\end{proposition}
The same argument used in Proposition~\ref{prop: d1=d2} can be used to prove the following claim.
\begin{proposition}
    \label{prop: cont-time chains}
    Let $F$ be a continuous-time semiflow on $X$ and set $f=F^1$.
    Then $\cC_{F,d}$, each $\cR_{\cC_{F,d,T}}$, $T>0$, and each $\cC_{f,d,N}$ for all $N>0$ are independent on the metric
    (which is why, in the items below, we omit the metric function in the indices of the chains relations).
    \begin{enumerate}          
        \item $G_F$ is a trapping region for $\cC_{F}$, for $\cC_{F,T}$ for all $T>0$ and for $\cC_{f,N}$ for all $N>0$;
        \item $\cR_{\cC_{F}}$, each $\cR_{\cC_{F,T}}$, $T>0$, and each $\cC_{f,N}$ for all $N>0$ are subsets of $G_F$;
        \item $\cR_{\cC_{F}}$, each $\cR_{\cC_{F,T}}$, $T>0$, their nodes and each $\cC_{f,N}$ for all $N>0$ are $F$-invariant.
    \end{enumerate}
\end{proposition}
In the remainder of the article, we will omit the metric function from the indices of the chains streams.

\black
Next two lemmas show that, for all that concerns infinitesimal chains, it is enough to consider the time-1 map $f=F^1$.
\begin{lemma}
    \label{lemma:long chains}
    For any $\eps>0$ and $x\in\cR_{\cC_{F,T}}$, there are $(F,\eps,T)$-chains of arbitrarily large length from $x$ to itself.
\end{lemma}
\begin{proof}
    Fix any integer $n>0$.
    By hypothesis, there is at least a $(F,\eps,T)$-chain $C$ from $x$ to itself.
    By concatenating $C$ with itself enough times, the result is a $(F,\eps,T)$-chain of length larger than $n$.
\end{proof}
In several statements below, starting from next one, we will use the notation ${\BF \lfloor T\rfloor}$ to indicate the largest integer not larger than $T$.
\begin{lemma}
    \label{lemma:TvsN}
    Let $x$ and $y$ be $\cC_{F,T}$-equivalent 
    and set $f=F^1$ and $N=\lfloor T\rfloor$.
    Then $x$ and $y$ are $\cC_{f,N}$-equivalent.
\end{lemma}
\begin{proof}
    By Proposition~\ref{prop: cont-time chains},
    $G_F$ is a trapping region for both streams $\cC_{F,T}$ and $\cC_{f,N}$.
    Hence, it is enough to consider the analogue problem for the restriction of $F$ to $G_F$.  
    Therefore, in the reminder of the proof we assume, without loss of generality, that $X$ is compact.
    
    We will prove that, for every $\eps>0$, there is a $(f,N,\eps)$-chain from $x$ to $y$. 
    The same argument then can be used to show that there is a $(f,\eps,N)$-chain from $y$ to $x$. 

    Fix an $\eps>0$ and let $\delta>0$ satisfy point (1) of Lemma~\ref{lemma: UC cont} and point (2) of Lemma~\ref{lemma: UC discrete} with $p=N$.
    We can assume without loss of generality that $\delta\leq\eps$.
    Let $C$ be a $(F,T,\delta)$-loop based at $x$ with points $(c_0,\dots,c_r)$ and times $(t_0,\dots,t_{r-1})$ such that 
    $c_i=y$ for some $0<i<r$.
    Set 
    $S=t_0+\dots+t_{r-1}$.
    Notice that, if $S/N$ is irrational, since rationals are dense, we can change $t_{r-1}$ to a new time $t'_{r-1}$ so that the new chain $(c_0,\dots,c'_{r-1},c_r)$ is still a $(F,T,\delta)$-loop based at $x$ but this time its period $S'$ is such that $S'/N$ is rational.
    Hence, we can assume without loss of generality that $S/N$ is rational.
    \black
    
    Following Hurley~\cite{Hur95}, we build a $(f,N,\eps)$-loop $C'$ based at $x$ 
    in the following way.
    For every $j$, set $s_j=\sum_{i=0}^jt_j$. 
    The $s_j$ are precisely the times at which, on the chain $C$, there are jumps -- precisely, a jump from $F^{t_j}(c_{j})$ to $c_{j+1}$.
    We start $C'$ by setting $c'_0=c_0$.
    Then, after Hurley, for each $k$, we follow the rule below:
    \begin{enumerate}
        \item if there is no jump in $C$ in the interval $(kN,(k+1)N]$, then we set $c'_{k+1}=f^N(c'_k)=F^N(c'_k)$;
        \item if there is a jump in $C$ at $s_j\in(kN,(k+1)N]$, then we set\\ $c'_{k+1}=F^{(k+1)N-s_j}(c_{j+1})$.
    \end{enumerate}
    The two cases above cover all possible cases because, since in $C$ jumps take place at least $T\geq N$ time units apart, there can be at most one jump in each interval $(kN,(k+1)N]$.

    Once the $(F,T,\delta)$-loop gets back to $x$, not necessarily this is the case for the $(f,N,\eps)$-chain.
    Nevertheless, recall that $S/N$ is rational, namely there are integers $m,n>0$ such that $mS=nN$.
    The chain obtained by repeating $m$ times the sequence $(c_0,\dots,c_r)$ and the relative times $(t_0,\dots,t_{r-1})$ is still a $(F,T,\delta)$-loop.
    Hence, after applying the construction above to this new loop, the next-to-last of the $(f,N,\eps)$-chain coincides with the next-to-last point of the $(F,T,\delta)$-loop.
    Since $\delta\leq\eps$, this shows that the $(f,N,\eps)$-chain can be completed to a $(f,N,\eps)$-look based at $x$ by adding $x$ as the last point of the chain.
\end{proof}
\begin{corollary}
    \label{cor:Ff}
    Assume that 
    $F$ is a continuous-time semi-flow with compact dynamics.
    Set $f=F^1$, fix a $T>0$ and set $N=\lfloor T\rfloor$.
    Then $\cR_{\cC_{F,T}}=\cR_{\cC_{f,N}}$ and each node of $\cR_{\cC_{F,T}}$ is a node of $\cR_{\cC_{f,N}}$ and viceversa.
\end{corollary}
\begin{lemma}
    Assume that 
    $F$ is a continuous-time semi-flow with compact dynamics.
    Set $f=F^1$, fix a $T>0$ and set $N=\lfloor T\rfloor$.
    Then, if $x$ is $\cC_{F,T}$-upstream of $y\in\cR_{\cC_{F,T}}$,
    $x$ is $\cC_{f,N}$-upstream of $y$.
\end{lemma}
\begin{proof}
    Let $C$ be a $(F,T,\eps)$-chain from $x$ to $y$ and let $D$ be a $(F,T,\eps)$-chain loop from $y$ to itself.
    Let $D^n$ be the concatenation of $D$ with itself $n$ times.
    Then by concatenating $C$ with $D^n$ we can get a chain from $x$ to $y$ of arbitrary length.
    Hence, by using the very same procedure of the previous lemma, we can prove that, for every $\eps>0$, there is a $(f,N,\eps)$-chain from $x$ to $y$.
\end{proof}
\begin{corollary}
    Assume that 
    $F$ is a continuous-time semi-flow with compact dynamics.
    Set $f=F^1$, fix a $T>0$ and set $N=\lfloor T\rfloor$.
    Then $\Gamma_{\cC_{F,T}}=\Gamma_{\cC_{f,N}}$.
\end{corollary}
The results above show already that all that the qualitative description of the dynamics of a continuous-time semi-flow $F$ with compact dynamics is all encoded in the powers of its time-1 map. 
Below we show that, in fact, the first power of the time-1 map is enough.

\begin{theorem}
    \label{thm: cont=disc}
    Let $f$ be a discrete-time semi-flow with compact dynamics.
    Then, for any integer $N>0$, $\cO_f\cup\cC_{f,N}=\cC_{f}$.
    If $f=F^1$ for some continuous-time semi-flow $F$, then we have also that $\cC_{F}=\cO_F\cup\cC_{f}$.
\end{theorem}
\begin{proof}
    As in the proof of Lemma~\ref{lemma:TvsN}, we can assume without loss of generality that $X$ is compact.

    First notice that, for any $N>0$, $\cC_{f}\subset\cO_f\cup\cC_{f,N}$ because every $(f,\eps,N)$-chain $C$ can be seen as a $(f,\eps)$-chain -- just break each jumpless segment in pieces of length 1.
    To complete the proof, we need to prove that, given any $N>0$ and $\eps>0$, if $x$ can be joined to $y$ by a $(f,\eta)$-chain for every $\eta>0$, then we can join $x$ to $y$ with a $(f,\eps,N)$-chain.

    So, let $\eps>0$, set $p=2N$ and let $\delta>0$ be the $\delta$ whose existence is granted by Lemma~\ref{lemma: UC discrete}(2) and $C$ a $(f,\delta)$-chain from $x$ to $y$.
    Recall that, by possibly concatenating $C$ with some $(f,\delta)$-chain from $y$ to itself, we can assume that $C$ has at least $N$ points.
    Let $c_0,\dots,c_r$, $r\geq N$, be the points of $C$.
    Then, by Lemma~\ref{lemma: UC discrete}(2), each pair $c_{kN},c_{(k+1)N}$ is a $(f,\eps,N)$-chain.
    If $r=q N$ for some integer $q>0$, then $c_0,c_N,\dots,c_{q N}$ is a $(f,\eps,N)$-chain from $x$ to $y$.
    Otherwise, $q N<r<(q+1)N$ for some $q>0$.
    In this case, we use as the final segment of the new chain the pair $c_{(q-1)N},c_r$. 
    Since $2N>r-(q-1)N>N$, even this pair is a $(f,\eps,N)$-chain.
    Hence, the chain $c_0,c_N,\dots,c_{(q-1)N},c_r$ is, in any case, a $(f,\eps,N)$-chain from $x$ to $y$.
    This proves that $\cO_f\cup\cC_{f,N}=\cC_{f}$.

    The second claim of the theorem comes from the fact that $\cC_{F,T} = \cO_F\cup\cC_{f,\lfloor T\rfloor}$ (Cor.~\ref{cor:Ff}) and that $\cO_F\cup\cC_{f,N}=\cO_F\cup\cC(f)$ for every integer $N>0$ (by the first claim of this theorem).
\end{proof}
\subsubsection{The {\BF$\Sigma$}chains streams}
\begin{definition}
    Given a discrete-time semiflow $F$ and a metric $d$ compatible with the topology of $X$, we call {\BF $(F,d)$-$\Sigma$chains stream} the relation
    \beq
    \Sigma_{F,d}  = \{(x,y):\text{for every $\eps>0$, there is a $(F,d,\eps)$-$\Sigma$chain from $x$ to $y$}\}.
    \eeq
\end{definition}
Unlike the chains streams, the $\Sigma$chains streams do depend on the metric even in case of compact dynamics, as the example below shows.
\begin{example}
    Let $X=[0,1]$ and let $F$ be a discrete-time flow on $X$ that has the ternary Cantor set $C$ as its set of fixed points and moves all other points rightwards, so that each of them asymptotes to the closest Cantor set point at its right.
    We claim that whether or not $1$ is $\Sigma_{F,d}$-downstream from $0$ depends on the metric.
    Indeed, the set of $\eps$-jumps of any $(F,d,\eps)$-chain from 0 to 1 must cover $C$, since points of $C$ are fixed, and so there is such a chain if and only if the length of $C$ with respect to $d$ is zero.
    In case of the Euclidean distance $d_E$, we know that the measure of $C$ is zero and so $(0,1)\in\Sigma_{F,d_E}$.
    Now, let $\varphi:[0,1]\to[0,1]$ be a homeomorphism such that the image of $C$ is the Smith-Volterra-Cantor set, which is a Cantor set of measure $1/2$.
    Then $d_\varphi(x,y)=d_E(\varphi(x),\varphi(y))$ is a metric on $[0,1]$ and, with respect to this metric, $C$ has measure $1/2$.
    Hence, $(0,1)\not\in\Sigma_{F,d_\varphi}$.
    
\end{example}
The following proposition is a direct consequence of the facts that $\Sigma_{F,d}\subset\cC_{F,d}$ and that every $(F,d,\eps)$-$\Sigma$chain is a $(F,d,\eps)$-chain.
\begin{proposition}
    \label{prop: Sigma_F}
    For a given semiflow $F$ with compact dynamics, denote by $F_G$ its restriction to $G_F$.
    The following hold for every metric $d$:
    \begin{enumerate}
        \item $\Sigma_{F,d}$ is a $\Omega$stream.
        \item $G_F$ is a trapping region for $\Sigma_{F,d}$.
        \item $\cR_{\Sigma_{F,d}}$ is $F$-invariant and all of its nodes are $F$-invariant.
        \item $\cR_{\Sigma_{F,d}}\subset G_F$.
        \item $\cR_{\Sigma_{F_G,d}}=\cR_{\Sigma_{F,d}}$;
        \item $N\subset\cR_{\Sigma_{F,d}}$ is a node of $\Sigma_{F,d}$ if and only if it is a node of $\Sigma_{F_G,d}$;
        \item $\Gamma_{\Sigma_{F_G,d}}=\Gamma_{\Sigma_{F,d}}$.
    \end{enumerate}
\end{proposition}
\subsubsection{The smallest stream}
\label{sec: A_F}
This stream was introduced by Joe Auslander in 1964~\cite{Aus63} as the smallest closed and transitive extension of the prolongational relation.
Below we prove that, under suitable assumptions, the smallest stream is a $\Sigma$chains stream.
\begin{proposition}
    Let $F$ be a semiflow with compact dynamics and
    denote by $\cM_X$ the set of all metrics on $X$ compatible with its topology.
    Then the following holds:
    \begin{enumerate}
        \item $\cA_F=\displaystyle\bigcap_{d\in\cM_X}\Sigma_{F,d}$;
        \item $\cR_{\cA_F}=\displaystyle\bigcap_{d\in\cM_X}\cR_{\Sigma_{F,d} }$;
        \item $\cA_F=\Sigma_{F,d} $ for some $d\in\cM_X$.
    \end{enumerate}
    In particular, $\cA_F$ satisfies all cases of Proposition~\ref{prop: Sigma_F}.
\end{proposition}
\begin{proof}
    Under the theorem's hypotheses, $G_F$ is a trapping region for all streams involved.
    Hence, outside of $G_F$ all these streams coincice and it is enough to prove that the properties in the claim hold within $G_F$.
    Since $G_F$ is compact, the proof of Proposition~6.4.1 in~\cite{DLY25} applies to it and so the claims follow.
\end{proof}
\subsubsection{Chains streams with countably many nodes}
\label{sec:countably many nodes}
\begin{definition}
    Assume $N_1,N_2$ are distinct nodes of a stream $S$.
    We say they are {\BF adjacent} if $N_1\Sto N_2$ and, whenever $N_1\Sto N\Sto N_2$, then either $N=N_1$ or $N=N_2$.
\end{definition}
\begin{lemma}
    \label{lemma:streams adjacent}
    Let $S$ be a substream of $\cC_F$.
    Then, there is a bitrajectory between every pair of $S$-adjacent nodes. 
\end{lemma}
\begin{proof}
    It is enough to prove the theorem for $S=\cC_F$.
    The same argument used in Proposition~\ref{prop: strong edge} applies to this case and shows that there is a bitrajectory $b$ with  $\alpha(b)\subset N_1$ and $\omega(b)\subset N_2$.
\end{proof} 
Next final result is proven by the same proof given in~\cite{DLY25}, except for the following update: for every $x\in X$, the reason why $\Omega_F(x)\neq\emptyset$ is that $F$ has a global attractor.
\begin{theorem}
    \label{thm: CF=AF}
    Let $F$ be a semiflow with compact dynamics.
    Then, if $\cA_F$ has countably many nodes, $\cA_F=\cC_F$. 
    In particular, $\Sigma_{F,d}=\cC_F$ for every metric $d$
    compatible with the topology of $X$.    
\end{theorem}

\section*{Acknowledgments}
The first author was partially supported by NSF grant \# 2308225.

\bibliographystyle{amsplain}  
\bibliography{refs}  
\end{document}